\documentclass[11pt, reqno]{article}
\usepackage{amsfonts}
\usepackage{amsmath}
\usepackage{amsthm}
\usepackage{amssymb}
\usepackage{mathrsfs}
\usepackage{graphicx}
\numberwithin{equation}{section}

\numberwithin{equation}{section}

\theoremstyle{plain}
\newtheorem{Theorem}{Theorem}[section]
\newtheorem{Lemma}[Theorem]{Lemma}
\newtheorem{Cor}[Theorem]{Corollary}

\theoremstyle{definition}
\newtheorem{Def}{Definition}[section]

\theoremstyle{remark}\newtheorem{Remark}{Remark}[section]

\newenvironment{Proof}{\noindent
{\bf\underline{Proof:} }}
{\hspace*{\fill}\qed\vskip1em}

\font\msbm=msbm10 \def\mathbb#1{\hbox{\msbm{#1}}}  
\newcommand{\N}{{\mathbb{N}}}
\newcommand{\R}{{\mathbb{R}}}
\newcommand{\Z}{{\mathbb{Z}}}
\newcommand{\C}{{\mathbb{C}}}
\newcommand{\TT}{{\mathbb{T}}}
\newcommand{\B}{{\mathcal{B}}}

\newcommand{\A}{{\mathcal{A}}}

\newcommand{\cJ}{{\mathcal{J}}}
\newcommand{\cS}{{\mathcal{S}}}

\renewcommand{\AA}{{\mathbb{A}}}
\newcommand{\BB}{{\mathbb{B}}}
\newcommand{\AwA}{\AA_w^\A}
\newcommand{\BwA}{\BB_w^\A}
\newcommand{\G}{{\mathcal{G}}}
\renewcommand{\L}{{\mathcal{L}}}

\newcommand{\K}{{\mathcal{K}}}
\renewcommand{\H}{{\mathcal{H}}}
\renewcommand{\Re}{{\operatorname{Re}}}
\newcommand{\HH}{{\mathbb{H}}}
\newcommand{\D}{{\mathcal{D}}}
\newcommand{\TPi}{{\widetilde{\pi}}}
\newcommand{\TV}{{\widetilde{V}}}
\newcommand{\STFT}{\operatorname{STFT}}

\newcommand{\Id}{\operatorname{Id}}
\newcommand{\Co}{{\mathsf{Co}}}

\newcommand{\rank}{{\operatorname{rank}}}

\newcommand{\ol}{\overline}

\newcommand{\spann}{\operatorname{span}}

\newcommand{\Kop}{{K}}

\newcommand{\bdiscrete}{\flat}
\newcommand{\ddiscrete}{\natural}

\begin{document}
\title{Radial Time-Frequency Analysis and Embeddings of Radial Modulation Spaces}
\author{Holger Rauhut\\University of Wroc{\l}aw, Mathematical Institute\\ 
Pl. Grunwaldzki 2/4, 50--384 Wroc{\l}aw, Poland\\
rauhut@ma.tum.de}

\maketitle

\begin{abstract}In this paper we construct frames of Gabor type for the space
$L^2_{rad}(\R^d)$ of radial $L^2$-functions, and more generally, for subspaces of modulation
spaces consisting of radial distributions.
Hereby, each frame element itself is a radial function. This construction is
based on a generalization of the so called Feichtinger-Gr\"ochenig theory -- sometimes
also called coorbit space theory -- 
which was developed in an earlier article. We show that this new type of 
Gabor frames behaves better in linear and non-linear approximation in a certain sense
than usual Gabor frames when approximating a radial function.
Moreover, we derive new embedding theorems for coorbit spaces restricted to invariant
vectors (functions) and apply them to modulation spaces
of radial distributions. As a special case this result implies that
the Feichtinger algebra $(S_0)_{rad}(\R^d) = M^1_{rad}(\R^d)$ restricted to radial
functions is embedded into the Sobolev space $H^{(d-1)/2}_{rad}(\R^d)$. Moreover,
for $d\geq 2$ the embedding $(S_0)_{rad}(\R^d) \hookrightarrow L^2_{rad}(\R^d)$ is compact.
\end{abstract}


\section{Introduction}

Nowadays time-frequency analysis is a well-developed 
field with many applications in signal analysis and wireless communication
\cite{GrBook}. The modulation spaces intruduced by Feichtinger in the early 80's
\cite{Fei_Mod} play a fundamental role in this mathematical area \cite{GrBook}.

The basic operators in time-frequency analysis on $\R^d$ 
are the translation $T_x f(y) = f(y-x)$ and the modulation 
$M_\omega f(y) = e^{2\pi i y\cdot \omega} f(y)$. 
For a fixed window function $g \in L^2(\R^d)$ the short time Fourier transform (STFT)
is given by $\STFT_g f(x,\omega) = \langle f, M_\omega T_x g\rangle$.
It is well-known that for
a suitable choice of constants $a,b >0$ and of $g \in L^2(\R^d)$ the family
$\{M_{bj} T_{ak} g:\, j,k \in \Z^d\}$ forms a frame - a so called Gabor frame.
In other words, the STFT admits discretizations.
In particular, we can write an arbitrary $L^2$-function as
$f = \sum_{j,k \in \Z^d} \lambda_{j,k} M_{bj} T_{ak} g$. Furthermore,
expansions of this type extend to modulation spaces \cite{GrBook}.

In \cite{Rau1,Rau2,Rau_Diss} we considered radial functions, or more generally functions which
are invariant under the action of some subgroup of $O(d)$, and raised the question whether
it is possible to develop an adapted time-frequency analysis for such functions. 
Of course, one can
apply all results that are valid for functions on $\R^d$. 
However, it seems natural that the additional information that the 
function under consideration is radial should allow one to gain some advantages.
When starting investigations in this direction one immediately observes that the operators
$T_x$ and $M_\omega$ do not preserve radiality except for trivial cases. 
So the natural requirement is 
to replace these operators by some that preserve radiality.
 
In \cite{Rau1} we found a natural candidate for such operators in the more general context of
square-integrable group representations. In particular,
we developed the Hilbert space theory for the
abstract continuous transform and applied it
to the special case of the short time Fourier transform (STFT) of radial 
functions.
In \cite{Rau2} we treated the corresponding discrete theory which lead to the construction of
discrete (Banach) frames for spaces of invariant vectors (functions). 
To do this we generalized the well-established coorbit space theory originally developed
by Feichtinger and Gr\"ochenig \cite{FG2,FG1,FG3,Gr}. This
discretization method does not only work on the Hilbert space level (i.e., for $L^2$)
but also for more general Banach spaces, so-called coorbit spaces of invariant vectors, 
see Section \ref{sec_notation}. 
In the special case of time-frequency analysis they are precisely the modulation spaces. 
We remark that coorbit space theory was further generalized to the setting of an 
abstract continuous frame in \cite{FR}.
   
In \cite{Rau2} we did not give examples. So in this paper we show in detail how to apply the
general results in order to obtain radial Gabor frames and atomic decompositions for radial
modulation spaces explicitly. In order to study the quality of approximation with radial Gabor
frames we compare it to the one with usual Gabor frames. In particular, we investigate linear and 
non-linear approximation. We first derive results in the general context of coorbit spaces 
and apply them to approximation in radial modulation spaces later on. 
To some extent it turns out that the radial Gabor frames really perform better than 
the usual Gabor frames
when approximating a radial function. This provides a further justification
of the original motivation to consider an adapted analysis for radial functions.

Moreover, we present how embeddings of 
coorbit spaces of invariant elements can be studied by means of embeddings of certain
sequence spaces. This result generalizes the work of Feichtinger and Gr\"ochenig in
\cite{FG2}.
The general embedding theorem indicates that
restricting coorbit spaces to invariant vectors (functions) may enforce certain embeddings
of coorbit spaces to become compact. In the case of the Besov and Triebel-Lizorkin spaces
this phenomenon was observed earlier by Skrzypczak et al. in \cite{KLSS,SickS,LS,ST}.

We apply these results to derive new embedding results for subspaces of modulation spaces
consisting of radial functions. Moreover, we determine the entropy numbers of certain
embeddings of radial modulation spaces.
As special case it turns out that
the Feichtinger algebra $(S_0)_{rad}(\R^d) = M^1_{rad}(\R^d)$ restricted to radial functions
is embedded
into the Sobolev space $H^{(d-1)/2}_{rad}(\R^d)$. This is rather surprising because $S_0$ (not restricted
to radial functions) is not a subspace of $H^{(d-1)/2}$. Moreover, for $d\geq 2$ 
the embedding $(S_0)_{rad}(\R^d) \hookrightarrow H^s_{rad}(\R^d)$ 
is compact if and only if $s<(d-1)/2$. In particular, $(S_0)_{rad}(\R^d)$ is compactly embedded into
$L^2_{rad}(\R^d)$ for $d\geq 2$.

The paper is organized as follows. In Section 2 we introduce the necessary notation
and background information. Moreover, we prove the general results about embeddings
of abstract coorbit spaces. Also we provide statements about linear and non-linear
approximation. In Section 3 we apply the general theory to time-frequency analysis
of radial functions in $\R^d$. We explicitly construct the radial Gabor frames
and state results about linear and non-linear approximation with these. Moreover,
we prove the mentioned embedding theorems for subspaces of modulation spaces
consisting of radial functions (distributions).

\section{Coorbit spaces of invariant elements}
\subsection{Notation and Preliminaries}
\label{sec_notation}

As announced we will derive our results first in the abstract
setting of coorbit spaces and then apply it to modulation spaces
of radial functions. To this end we need to introduce some notation
from \cite{FG1,Rau1,Rau2,Rau_Diss}.

Let $\G$ be a locally compact group and $\A$ be a compact 
automorphism group of $\G$, such that $\A$ acts continuously 
on $\G$, i.e., the mapping 
$\G \times\A  \to \G, (x,A) \mapsto Ax$ is continuous. 
We denote the left Haar measures on $\G$ and $\A$ by
$\mu$ and $\nu$, where $\nu$ is assumed to be normalized. 
However, we usually write $dx$ and $dA$ in integrals. 
The modular function on $\G$ is denoted 
by $\Delta$ and the left and right translation operators on $\G$ 
by $L_y F(x) = F(y^{-1} x)$ and $R_y F(x) = F(xy)$. Furthermore, we define
two involutions by $F^\vee(x) = F(x^{-1})$ and $F^\nabla(x) = \ol{F(x^{-1})}$. 
The action of $\A$ on functions on $\G$ is denoted by $F_A(x) = F(A^{-1}x)$, $A \in \A$, and
the action on measures $\tau \in M(\G)$, the space of complex bounded Radon measures on $\G$ 
(the dual space of $C_0(\G)$), by $\tau_A(F) = \tau(F_{A^{-1}})$, $A \in \A$, $\tau \in M(\G)$, $F\in C_0(\G)$.

The functions (measures) which satisfy $F_A = F$ for all $A \in \A$ are called invariant (under $\A$).
A standard argument shows that the Haar-measure $\mu$ and the modular function $\Delta$ are invariant
under any compact automorphism group \cite{Rau_Diss}. 
For a function (measure) space $Y$ on $\G$ we denote its subspace of 
invariant elements by $Y_\A:=\{F \in Y, F_A=F \mbox{ for all } A \in \A\}$. 
An invariant function on $\G$ can be interpreted
as a function on $\K:= \A(\G)$ the space of all orbits of the 
form $\A x:= \{Ax, A \in \A\}, x\in \G$. The orbit space $\K$ becomes 
a hypergroup by 
inheriting the topology  and the convolution structure 
of $\G$ in a natural way \cite{Jew,Rau1,Rau_Diss}. 

For some positive measurable weight function $m$ on $\G$ we define the weighted space
$L^p_m:=\{F \mbox{ measurable}, Fm \in L^p\}$ with norm $\|F|L^p_m\|:= \|Fm|L^p\|$
where the $L^p$-spaces on $\G$ are defined as usual. 

In this paper we will work with Banach spaces of functions 
on $\G$ which will usually be denoted by $Y$. Similarly as in \cite{Gr} we will 
make the following assumptions on $Y$.
\begin{enumerate}
\item $Y$ is continuously embedded into $L^1_{loc}(\G)$, the locally integrable functions on $\G$.
\item $Y$ is solid, 
i.e., if $F \in L^1_{loc}(\G), G \in Y$ and $|F(x)| \leq |G(x)|$ a.e. then
$F \in Y$ and $\|F|Y\| \leq \|G|Y\|$.
\item  $Y$ is invariant under left and right translations. We may hence define the two functions
$u(x) := \|L_x|Y\to Y\|$ and $v(x):= \|R_{x^{-1}}|Y \to Y\| \Delta(x^{-1})$. Clearly, 
$u(xy) \leq u(x)u(y)$ and $v(xy) \leq v(x) v(y)$, i.e., $u$ and $v$ are submultiplicative.
Additionally, we require that $u$ and $v$ are continuous. Under these assumptions it holds,
see \cite{FG1,Rau_Diss}
\begin{equation}\label{conv_r}
Y *L^1_v \subset Y,\quad \|F*G|Y\| \leq \|F|Y\| \,\|G|L^1_v\|
~\mbox{ for all }F \in Y, G\in L^1_v.
\end{equation}
\item $\A$ acts continuously on $Y$. Without loss of generality we may 
then even assume that $\A$ acts isometrically on $Y$ implying
 $u(Ax) = u(x)$ and $v(Ax) = v(x)$
for all $A\in \A$. (In case this is not true define an invariant norm on $Y$ by
$\|F|Y\|' := \int_\A \|F_A|Y\| dA$. Since $\A$ acts 
continuously on $Y$ this is an equivalent norm on $Y$.)
Then $Y_\A$ is a closed non-trivial subspace of $Y$. 
\end{enumerate}
Given a submultiplicative weight $w$, another continuous 
weight function $m$ is called
$w$-moderate if $m(xyz) \leq w(x) m(y) w(z)$ for all $x,y,z \in \G$.
Spaces $L^p_m$ with invariant moderate weight function $m$ are examples 
of spaces $Y$ with the properties above.

We will always associate a weight function $w$ to $Y$ which is defined by
\begin{equation}\label{def_weight}
w(x) \,:=\, \max\{u(x), u(x^{-1}), v(x), v(x^{-1})\Delta(x^{-1})\}.
\end{equation}
Then as a consequence $w$ is continuous, 
$w(xy)\leq w(x) w(y)$, $w(x) \geq 1$ and $w(Ax) = w(x)$ for all $A\in \A$ and $x\in \G$.
Furthermore, by (\ref{conv_r}) it holds 
\begin{equation}\label{Yw_conv}
Y * L^1_w \,\subset\, Y, \quad \|F*G|Y\| \,\leq\, \|F|Y\|\,\|G|L^1_w\|.
\end{equation}

We further assume that we have given a unitary, irreducible (strongly continuous) representation $\pi$ of
$\G$ on some Hilbert space $\H$ and some unitary (strongly continuous) representation $\sigma$ of $\A$
on the same Hilbert space $\H$ (not necessarily irreducible) 
such that the following basic relation is satisfied 
(see also \cite{Rau1,Rau_Diss}),
\begin{equation}\label{Cond_pi_sigma}
\pi(A(x)) \sigma(A) \,=\, \sigma(A) \pi(x) \qquad\mbox{ for all }x\in \G, A\in \A.
\end{equation}
In other words, we require that the representations $\pi_A:= \pi \circ A$
are all unitarily equivalent to $\pi$ and that the intertwining operators
$\sigma(A)$ form a representation of $\A$.

For $f \in \H$ we denote $f_A = \sigma(A) f$ and $\H_\A := \{f \in \H, f_A = f \mbox{ for all } A \in \A \}$,
the closed(!) subspace of invariant elements.
We always assume that $\H_\A$ is not trivial.
The wavelet transform or voice transform is defined by
\[
V_g f (x) \,:=\, \langle f, \pi(x) g\rangle.
\]
It maps $\H$ into $C^b(\G)$, the space of bounded continuous functions on $\G$. 
With an element $g\in \H_\A$ we denote by $\tilde{V}_g$ the restriction of $V_g$ to $\H_\A$. 
We recall some facts from \cite{Rau1,Rau_Diss}.
\begin{itemize}
\item For $f,g \in \H_\A$ the function $\TV_g f$ is invariant under 
$\A$, i.e., $\TV_g$ maps $\H_\A$ into $C^b_\A(\G)$.
\item For $x \in \G$ we define  
\begin{align}\label{def_TPi}
\TPi(x) \,:=\, \int_\A \pi(A x) dA
\end{align}
in a weak sense. This operator maps $\H_\A$ onto $\H_\A$ and 
depends only on the orbit of $x$ under $\A$, i.e.,
$\TPi(B x) = \TPi(x)$ for all $B \in \A$. Furthermore, it holds 
\begin{equation}\label{wave_inv_comp}
\TV_g f(x) \,=\, \langle f, \TPi(x) g\rangle_{\H_\A} .
\end{equation}
\item The operators $\TPi(x)$ form an irreducible representation of the orbit hypergroup $\K$.
\end{itemize}
We further require that $\pi$ is integrable which means that there exists a nonzero element $g \in \H$
such that $\int_\G|V_g g(x)|dx < \infty$. This 
implies that $\pi$ is square-integrable, i.e., 
there exists $g \in \H$ such that
$\int_\G |V_g f(x)|^2 dx < \infty$ for all $f \in \H$. Such a $g$ (corresponding to the
square-integrability condition) is called admissible. 
We list some further properties from \cite{DM} and \cite{Rau1} that hold under 
the square-integrability condition.
\begin{itemize}
\item There exists a positive, densely defined operator $\Kop$ such that the domain $\D(\Kop)$ of $\Kop$ consists
of all admissible vectors and the orthogonality relation
\[
\int_\G V_{g_1} f_1(x) \ol{V_{g_2} f_2(x)} dx \,=\, \langle \Kop g_2, \Kop g_1 \rangle \langle f_1, f_2 \rangle
\]
holds for all $f_1,f_2 \in \H, g_1, g_2 \in \D(\Kop)$.
\item The operator $\Kop$ commutes with the action of $\A$, i.e., 
$\sigma(A) \Kop = \Kop \sigma(A)$ for all $A \in \A$. Furthermore, $\D_\A(\Kop):= \D(\Kop) \cap \H_\A$ is dense
in $\H_\A$ and $\Kop$ maps $\D_\A(\Kop)$ into $\H_\A$.
\item For $g\in \D_\A(\Kop)$ with $\|\Kop g\|=1$ we have the following inversion formula on $\H_\A$
\begin{equation}\label{Vg_inv}
f \,=\, \int_\K \TV_g f(y) \TPi(y) g \,d\widetilde{\mu}(y),\quad f\in \H_\A
\end{equation}
where $\widetilde{\mu}$ denotes the projection of the Haar measure $\mu$ 
onto the orbit space $\K$.
The integral is understood in a weak sense.
\end{itemize}

In order to define the coorbit spaces we need to introduce the space of analyzing vectors.
For some submultiplicative weight function $w$ satisfying $w\geq 1$ it is defined by
\[
\AA_w \,:=\, \{ g \in \H,\,V_g g \in L^1_w(\G) \}
\]
and its subspace of invariant elements by
$\AwA \,:=\, \AA_w \cap \H_\A$. 
We only consider those weights for which $\AwA \neq \{0\}$. It is clear that
$\AwA \subset \D(\Kop)$. Now for some fixed non-zero vector $g \in \AwA$ we define
\[
\H^1_w \,:=\, \{f \in \H,\,V_g f  \in L^1_w \}
\]  
with norm 
\[
\| f|\H^1_w\| \,:=\, \|V_g f | L^1_w\|.
\] 
Its subspace of invariant elements is denoted by
\[\label{def_H1wA}
(\H_w^1)_\A \,:=\, \H_\A \cap \H^1_w \,=\, \{f \in \H_\A,\,\TV_g f  \in L^1_w\}.
\]
It is shown in \cite{Rau2,Rau_Diss} that $\TPi(x) g \in \AwA$ for all $g \in \AwA, x\in \G$. 
As a consequence, 
both $\AwA$ and $(\H^1_w)_\A$ are dense in $\H_\A$. As a reservoir for the general
coorbit spaces we take the space $(\H^1_w)^\urcorner$
of all bounded conjugate linear functionals on $\H^1_w$ (the anti-dual of $\H^1_w$) 
and its subspace $(\H^1_w)_\A^\urcorner$ of invariant elements (the anti-dual of 
$(\H^1_w)_\A$), respectively.
We may extend the voice transform onto $(\H^1_w)^\urcorner$ 
by
\[
V_g f(x) \,=\, f(\TPi(x) g) \,=\, \langle f, \TPi(x) g\rangle, \quad f \in \H^1_w, g \in \AA_w.
\]

Now let $Y$ be some function space on $\G$ that satisfies our hypothesis and let
$w$ be the weight function defined by (\ref{def_weight}). Then for some
fixed non-zero vector $g \in \AwA$ the coorbit space associated to $Y$ is defined by
\[
\Co Y \,:=\, \{ f\in (\H^1_w)^\urcorner, V_g f \in Y \}
\]
with natural norm
\[
\|f|\Co Y\| \,:=\, \|V_g f| Y\|.
\]
Its subspace of invariant elements is denoted by
\[
\Co Y_\A \,:=\, \Co Y \cap (\H^1_w)_\A^\urcorner\,=\, \{f \in (\H^1_w)_\A^\urcorner, 
\TV_g f \in Y_\A\}.
\] 
It was shown in \cite{FG1,Rau2,Rau_Diss} that the coorbit spaces are Banach spaces
whose definition does not depend on the particular choice of $g \in \AwA$.

We will write $A \asymp B$ if there exist constants $C_1,C_2>0$ such that 
$C_1 A \leq B \leq C_2 B$ independently of other expressions on which $A,B$ might depend.

\subsection{Atomic decompositions and Banach frames in coorbit spaces of invariant elements}

In \cite{Rau2,Rau_Diss} atomic decompositions and Banach frames for 
coorbit spaces of invariant elements have been derived. 

In order to state these results 
we need to recall a definition from \cite{Rau2,Rau_Diss}.

\begin{Def} Let $X=(x_i)_{i\in I}\subset \G$ be some family of points in $\G$ indexed
by some discrete index set $I$ and let $V=V^{-1}=\A(V)$ be some relatively compact
neighborhood of $e \in \G$.
\begin{itemize}
\item[(a)] $X$ is called $V$-$\A$-dense  if $\G = \bigcup_{i\in I} \A(x_i V)$.
\item[(b)] $X$ is called $V$-$\A$-separated if $\A(x_i V) \cap \A(x_j V) = \emptyset$ for all $i \neq j$.
\item[(c)] $X$ is called relatively separated with respect to $\A$ if for all compact sets
$W \subset \G$ there exists some constant $C_W$ such that
\[
\sup_{j \in I} \#\{ i\in I,\, \A(x_iW) \cap \A(x_jW) \neq \emptyset \} \,\leq\, C_W \,<\, \infty.
\] 
\item[(d)] $X$ is called well-spread with respect to $\A$ if it is both $V$-$A$-dense (for some $V$)
and relatively separated.
\end{itemize}
\end{Def}
The existence of well-spread sets was shown in \cite{Rau_Diss,Rau1}.
Moreover, we have the following lemma 
relating separated and relatively separated families.
\begin{Lemma}\label{lem_partition} 
Let $X=(x_i)_{i \in I} \subset \G$. The following properties
are equivalent.
\begin{itemize}
\item $X=(x_i)_{i\in I}$ is relatively separated with respect to $\A$.
\item For any compact set $K = K^{-1} = \A(K)$ there exists a finite partition
of the index set $I = \bigcup_{r = 1}^s I_r$ such that each family
$(x_i)_{i\in I_r}$ is $K$-$\A$-separated.
\end{itemize}
\end{Lemma}
\begin{Proof} This follows
immediately from Lemma 2.9 in \cite{FGr}.
\end{Proof}

Let us now define the sequence spaces that will characterize the coorbit spaces
$\Co Y_\A$. For some well-spread set $X=(x_i)_{i\in I}$ let
\begin{align}
Y_\A^\bdiscrete \,:=\, Y_\A^\bdiscrete(X) \,:=&\, \{ (\lambda_i)_{i \in I}, \sum_{i\in I} 
|\lambda_i|\chi_{\A(x_i U)} \in Y\},
\notag\\
Y_\A^\ddiscrete \,:=\, Y_\A^\ddiscrete(X) \,:=&\,\{ (\lambda_i)_{i \in I}, \sum_{i\in I} |\lambda_i| 
|\A(x_i U)|^{-1} \chi_{\A(x_i U)} \in Y\}\notag
\end{align}
with natural norms
\begin{align}
\|(\lambda_i)_{i \in I} | Y_\A^\bdiscrete\| \,:=&\, \|\sum_{i \in I} |\lambda_i| \chi_{\A(x_i U)}|Y\|,\notag\\
\|(\lambda_i)_{i \in I} | Y_\A^\ddiscrete\| \,:=&\, \|\sum_{i \in I} |\lambda_i||\A(x_i U)|^{-1} \chi_{\A(x_i U)}|Y\|.
\notag
\end{align}
Hereby, $|\A(x_i U)|$ denotes the Haar measure 
and $\chi_{\A(x_i U)}$ the characteristic function 
of the set $\A(x_i U)$. If $Y= L^p_m(\G), 1\leq p \leq \infty$, 
with invariant moderate weight function 
$m$ then $Y_\A^\bdiscrete(X) = \ell^p_{m_p}(I)$ 
and $Y_\A^\ddiscrete(X) = \ell^p_{\nu_p}(I)$ (with equivalent norms) where
\[
m_p(i) \,:=\, m(x_i) |\A(x_i U)|^{\frac{1}{p}}, \qquad \nu_p(i) \,:=\, m(x_i) |\A(x_i U)|^{\frac{1}{p}-1} .
\] 

As another ingredient we need Wiener amalgam spaces as introduced by
Feichtinger \cite{FeiW1,FeiW2}. For the definition we take a two-sided translation invariant
solid BF-space $Y$ and another two-sided invariant Banach space $B$ of functions or measures
on $\G$. Using a non-zero window function $k \in C_c(\G)$ we define the control function
\begin{equation}\label{def_control}
K(F,k,B)(x) \,:=\,\|(L_x k) F|B\|, \quad x \in \G
\end{equation}
where $F$ is locally contained in $B$, in symbols $F \in B_{loc}$. 
The Wiener amalgam space $W(B,Y)$ is then defined by
\[\label{def_Wiener_space}
W(B,Y) \,:=\, \{F \in B_{loc},\, K(F,k,B) \in Y\}
\]
with norm 
\[
\|F|W(B,Y)\|\,:=\, \|K(F,k,B)|Y\|.
\]
It was shown in \cite{FeiW2} that these spaces are two-sided invariant
Banach spaces which do not depend on the particular choice of the window
function $k$. Moreover, different functions $k$ define equivalent
norms. We will mainly need the spaces $W(L^\infty,Y)$ and $W(C_0,Y)$.
Replacing the left translation $L_x$ with the right translation $R_x$
in the definition (\ref{def_control}) of the control function leads to right Wiener
amalgam spaces $W^R(B,Y)$. 
The subspace of functions which are invariant under $\A$ is denoted
by $W_\A(B,Y)$ (or $W_\A^R(B,Y)$ for right amalgams). If $\A$ acts
isometrically on $B$ and $Y$ then $W_\A$ is a closed subspace
of $W_\A(B,Y)$.

In order to derive Banach frames we need the 'better space' of analyzing vectors
\[
\BB_w := \{g \in \AA_w, V_g g \in W^R(C_0,L^1_w)\},\qquad \BwA := \{g \in \AwA, \TV_g g \in W^R_\A(C_0, L^1_w)\}. 
\] 
It was shown in \cite{FG2,Rau_Diss} that $\BB_w$ is dense in $\H$ and that $\BwA$ is dense in $\H_\A$.

Now we are ready to formulate the result
concerning Banach frames and atomic decompositions of $\Co Y_\A$, see
\cite[Theorem 7.3]{Rau2} or \cite[Theorem 4.6.3]{Rau_Diss}.

\begin{Theorem}\label{thm_ab} Let $g \in \BwA \setminus \{0\}$. 
Then there exists a relatively compact neighborhood $U=U^{-1}=\A U$ of $e \in \G$ such that 
for any family 
$X=(x_i)_{i\in I}$, which is $U$-dense and well-spread with respect to $\A$, 
the family $\{\TPi(x_i) g\}_{ i\in I}$ is an atomic decomposition
for $\Co Y_\A$. This means that
\begin{itemize} 
\item there exist elements $\{e_i\}_{i \in I}$ in $(\H^1_w)_\A$ such 
that $(\langle f, e_i\rangle )_{i \in I} \in Y^\ddiscrete_\A(X)$ for all 
$f \in \Co Y_\A$ and 
\[
\|(\langle f, e_i\rangle)_{i\in I} | Y_\A^\ddiscrete\| \,\asymp\, \|f|\Co Y_\A\|;
\]   
\item it holds $f=\sum_{i\in I} \langle f, e_i \rangle \TPi(x_i) g$ for all 
$f \in \Co Y_\A$ with norm convergence if the finite sequences are dense
in $Y^\ddiscrete_\A$ and with w-$*$ convergence in $(\H^1_w)_\A^\urcorner$
in general.
\end{itemize}
Moreover, $\{\TPi(x_i)g\}_{i\in I}$ is a Banach frame for $\Co Y_\A$, i.e.,
\begin{itemize}
\item $(\langle f, \TPi(x_i)g\rangle)_{i\in I} \in Y^\bdiscrete_\A$ for all 
$f \in \Co Y_\A$ and 
\[
\|(\langle f, \TPi(x_i)g\rangle )_{i\in I}|Y^\bdiscrete_\A\| 
\, \asymp \, \|f|\Co Y_\A\|;
\]
\item there exists a bounded operator $\Omega : Y^\bdiscrete_\A \to \Co Y_\A$ 
such that $\Omega (\langle f, \TPi(x_i)g \rangle)_{i\in I} = f$ for all 
$f \in \Co Y_\A$.  
\end{itemize}
\end{Theorem}
This is a very general discretization theorem. Moreover, it characterizes
invariant coorbit spaces by means of the sequence spaces $Y^\ddiscrete_\A$
and $Y^\bdiscrete_\A$. We remark that in \cite{Rau2,Rau_Diss} we 
stated an explicit condition on the
set $U$ which uses a certain maximal function. This condition depends only on $g$ and on
$w$. This means that the theorem is valid 'uniformly' for all $Y$ whose associated
weight function $w_Y$ defined in (\ref{def_weight}) is dominated by $w$.

We will apply Theorem \ref{thm_ab} to the
special case of radial Gabor frames in the next section. For an application to
radial wavelet frames we refer to \cite{Rau_Diss}.

\subsection{Embeddings of Coorbit Spaces}

Next we study how embeddings of coorbit spaces
and embeddings of sequence spaces are related.
We start with an auxiliary result about interpolation
of voice transforms.

\begin{Theorem}\label{thm_interpol} 
Suppose $Y$ and $w$ are related as usual and let $g \in \BwA$ with
$\|Kg|\H\|=1$. Then there exists some compact set $K$ 
such that for any $K$-$\A$-separated 
family $X=(x_i)_{i \in I_r}$ there exists a linear bounded 
operator $S: Y^\bdiscrete_\A(X) \to \Co Y_\A$ such that whenever
$f = S (\lambda_i)_{i\in I}$ for 
$(\lambda_i)_{i\in I} \in Y^\bdiscrete_\A(X)$ then 
\[
(V_g f(x_i))_{i \in I} = (\lambda_i)_{i \in I} \quad \mbox{and}\quad
\|f|\Co Y\| \leq C \|(\lambda_i)_{i \in I}|Y^\bdiscrete_\A(X)\|.
\]
\end{Theorem}
\begin{Proof} This Theorem is shown in completely the same way as 
Proposition 8.2 in \cite{FG3}.
In particular, this task consists in slightly adjusting the proof
of Theorem 7.3 in \cite{FG3}. We omit the details.
\end{Proof}

In other words the mapping $\Co Y_\A \to Y^\bdiscrete_\A(X)$,
$f \mapsto (V_g f (x_i))_{i\in I}$
is surjective if $X$ is $K$-$\A$-separated for $K$ large enough.
Further, we need
some statements about
Wiener amalgam spaces of invariant elements.

\begin{Lemma}\label{lemA} Let $X=(x_i)$ be some well-spread set with respect to
$\G$ and $Q=Q^{-1}=\A(Q)$ be some neighborhood of $e \in \G$. 
Then a function 
$F$, which is invariant under $\A$, is contained
in $W_\A(L^\infty,Y)$ if and only if 
$(\|\chi_{\A(x_i Q)} F\|_\infty)_{i \in I} \in Y^\bdiscrete_\A(X)$ 
and there are constants
$C_1,C_2 > 0$ such that
\begin{equation}\label{WA_equiv}
C_1 \|F|W_\A(L^\infty,Y)\| \,\leq\, 
\|(\|\chi_{\A(x_i Q)} F\|_\infty)_{i \in I}|Y^\bdiscrete_\A(X)\|  
\,\leq\, C_2 \|F|W_\A(L^\infty,Y)\|.
\end{equation}
\end{Lemma}
\begin{Proof} Let us assume that the characteristic function 
$\chi_Q$ is taken for the definition of $W_\A(L^\infty,Y)$.

Suppose $F \in W_\A(L^\infty,Y)$. We note that
\[
\sum_{i \in I} \|F \chi_{\A(x_i Q)}\|_\infty \chi_{\A(x_i Q)}(x)
\,=\, \sum_{i\in I_x} \sup_{y \in \A(x_i Q)}|F(y)| \chi_{\A(x_i Q)}(x)
\]
where the sum runs over the finite index set $I_x = \{i \in I, x \in \A(x_i Q)\}$.
Indeed, $\# I_x \leq N$ uniformly in $x$. This yields
\[
\sum_{i \in I} \|F \chi_{\A(x_iQ)}\|_\infty \chi_{\A(x_i Q)}(x)
\leq N \sup_{y \in \A(x Q^2)}|F(y)|.
\]
By invariance of $F$ under $\A$ it holds
\[
K(F,\chi_{Q^2},L^\infty)(x) \,=\, \sup_{y \in xQ^2}|F(y)| \,=\, 
\sup_{y \in \A(xQ^2)}|F(y)|.
\]
Using the solidity of $Y$ we obtain 
\[
\|(\|F\chi_{\A(x_i Q)}\|_\infty)_{i\in I}|Y^\bdiscrete_\A\| 
\,\leq\,N \|K(F,\chi_{Q^2},L^\infty)|Y\| \leq C_1 \|F|W_\A(L^\infty,Y)\|.
\]   
Hereby, it is used that different window functions generate equivalent
norms on $W(L^\infty,Y)$.

For the converse inequality, let 
\[
J_x := \{i \in I, \A(x_iQ) \cap \A(x Q) \neq \emptyset\}.
\]
By the defining properties of a well-spread set $J_x$ is finite
and $\# J_x \leq N$ uniformly in $x$. Moreover, since the sets 
$\A(x_iQ)$ form a covering of $\G$ it holds
$\A(xQ) \subset \bigcup_{i \in J_x} \A(x_i Q)$.
By invariance of $F$ under $\A$ we obtain
\begin{align}
K(F,\chi_Q,L^\infty)(x) \,&=\, \sup_{y \in xQ}|F(y)| 
\,=\, \sup_{y \in \A(xQ)}|F(y)| 
\,\leq\, \sum_{i \in J_x} \sup_{y \in \A(x_i Q)} |F(y)|\notag\\
&\leq\, \sum_{i \in I} \sup_{y \in \A(x_i Q)}|F(y)| \chi_{\A(x_i Q^2)}(x).\notag 
\end{align}
The last inequality follows from the fact that 
$\A(xQ) \cap \A(x_iQ) \neq \emptyset$ implies $x \in \A(x_i Q^2)$.
By solidity of $Y$ this gives the lower estimate in (\ref{WA_equiv}).
\end{Proof}

\begin{Lemma}\label{lem_embed_equiv} 
Suppose $Y,Z$ are two solid BF spaces on $\G$. Then
$W_\A(L^\infty,Y) \subset W_\A(L^\infty,Z)$
if and only if $Y^\bdiscrete_\A(X) \subset Z^\bdiscrete_\A(X)$.
\end{Lemma} 
\begin{Proof} First assume $Y^\bdiscrete_\A(X) \subset Z^\bdiscrete_\A(X)$ 
and let $F \in W_\A(L^\infty,Y)$. By Lemma \ref{lemA} 
$(\|F\chi_{\A(x_iQ)}\|_\infty)_{i \in I} \in Y^\bdiscrete_\A(X) \subset
Z^\bdiscrete_\A(X)$ and, thus, again by Lemma \ref{lemA} $F$ is 
contained in $W_\A(L^\infty,Z)$. 

Conversely, suppose $W_\A(L^\infty,Y) \subset W_\A(L^\infty,Z)$ 
and let
$(\lambda_i)_{i \in I} \in Y^\bdiscrete_\A(X)$. We claim that the function
$F(x) = \sum_{i \in I} |\lambda_i| \chi_{\A(x_i Q)}(x)$ is contained
in $W_\A(L^\infty,Y)$. Indeed,
\begin{align}
K(F,\chi_Q,L^\infty)(x) \,&=\, 
\sup_{y \in xQ} \sum_{i \in I} |\lambda_i| \chi_{\A(x_iQ)}(y)
\,\leq\, \sum_{i \in I} |\lambda_i| \sup_{y \in xQ} \chi_{\A(x_iQ)}(y)\notag\\
&\leq\, \sum_{i \in I} |\lambda_i| \chi_{\A(x_iQ^2)}(x).\notag
\end{align}
The latter function is contained in $Y$ by definition of $Y^\bdiscrete_\A(X)$,
hence $F \in W_\A(L^\infty,Y)$. Since $W_\A(L^\infty,Z) \subset Z_\A$
we thus deduce $F \in Z_\A$ which is 
equivalent to $(\lambda_i)_{i\in I} \in Z^\bdiscrete_\A(X)$ by
definition. Altogether we deduced 
$Y^\bdiscrete_\A(X) \subset Z^\bdiscrete_\A(X)$. 
\end{Proof}

\begin{Cor}\label{Cor_unamb} 
Suppose $X_1$ and $X_2$ are two well-spread sets and let $Y,Z$
be solid BF-spaces. Then
$Y^\bdiscrete_\A(X_1) \subset Z^\bdiscrete_\A(X_1)$ if and only if
$Y^\bdiscrete_\A(X_2) \subset Z^\bdiscrete_\A(X_2)$. Hence,
we may unambiguously write $Y^\bdiscrete_\A \subset Z^\bdiscrete_\A$.
\end{Cor}
\begin{Proof} The assertion follows easily from
Lemma \ref{lem_embed_equiv} since $W_\A(L^\infty,Y)$ does not depend
on the well-spread set $X$.
\end{Proof}

Now we are ready to prove the following characterization of coorbit spaces.

\begin{Theorem}\label{thm_seq_coorbit} It holds $\Co Y_\A \subset \Co Z_\A$ if and only if 
$Y^\bdiscrete_\A \subset Z^\bdiscrete_\A$. 
In particular, two coorbit spaces
coincide if and only if the corresponding sequence spaces coincide.
\end{Theorem}
\begin{Proof} Assume $Y_\A^\bdiscrete \subset Z_\A^\bdiscrete$. By Lemma
\ref{lem_embed_equiv} this implies $W_\A(L^\infty,Y) \subset W_\A(L^\infty,Z)$.
It follows from Theorem 8.3 in \cite{FG2} that $\Co Y_\A = \Co W_\A(L^\infty, Y)$. Since
the implication $Y^1_\A \subset Y^2_\A \Longrightarrow \Co Y^1_\A \subset \Co Y^2_\A$ 
is trivial we conclude that $\Co Y_\A \subset \Co Z_\A$.

Conversely, let $\Co Y_\A \subset \Co Z_\A$. By Corollary \ref{Cor_unamb} 
we may choose the well-spread set $X=(x_i)_{i\in I}$ according to Theorem \ref{thm_interpol}.
Now assume $Y_\A^\bdiscrete(X) \not\subset Z_\A^\bdiscrete(X)$, i.e., that there
exists some $(\lambda_i)_{i\in I} \in Y_\A^\bdiscrete(X) \setminus Z_\A^\bdiscrete(X)$.
We have $f = S (\lambda_i)_{i \in I} \in \Co Y_\A \subset \Co Z_\A$ with $S$ as in
Theorem \ref{thm_interpol}. By \cite[Theorem 5.8]{Rau2} (see also 
\cite[Theorem 4.5.15]{Rau_Diss}) we conclude
$(\lambda_i)_{i\in I} = (V_g f (x_i))_{i\in I} \in Z_\A^\bdiscrete(X)$, a
contradiction. Thus, $Y^\bdiscrete_\A \subset Z^\bdiscrete_\A$. 
\end{Proof} 

We remark that the embeddings in the previous theorems are automatically continuous
by Theorem \ref{thm_coorbit_embed} below. 

Now we are ready to show the double retract property of the coorbit spaces
of invariant elements. Recall that some Banach space $B_1$ is called
a retract of the Banach space $B_2$ if there exist bounded linear 
operators $S: B_1 \to B_2$ and $T:B_2 \to B_1$ such that $T\circ S = \Id_{B_1}$.

\begin{Theorem} Choose $X=(x_i)_{i \in I}$ as in Theorem \ref{thm_ab}, 
i.e., such that $\{\TPi(x_i)g\}_{i \in I}$ is a Banach frame for $\Co Y_\A$.
Then $\Co Y_\A$ is a retract of $Y^\bdiscrete_\A(X)$ and, conversely,
$Y_\A^\bdiscrete(X)$ is a retract of a finite direct sum of copies
of $\Co Y_\A$.
\end{Theorem}
\begin{Proof} Denote by $A: \Co Y_\A \to Y^\bdiscrete_\A$ the operator
$f \mapsto (\langle f, \TPi(x_i) g\rangle)_{i \in I}$. Since $\{\TPi(x_i)g\}_{i \in I}$
is a Banach frame with 
bounded reconstruction operator $\Omega: Y^\bdiscrete_\A \to \Co Y$ we have
$\Omega \circ A = \Id_{\Co Y_\A}$, i.e., $\Co Y_\A$ is a retract of 
$Y^\bdiscrete_\A$. 

For the converse we choose a compact set $K=K^{-1}=\A(K)$ according to 
Theorem \ref{thm_interpol} and split the index set $I$ into finitely many subfamilies
$I_r$, $r=1,\hdots,s$, such that $(x_i)_{i \in I_r}$ is $K$-$\A$-separated for all $r$,
see Lemma~\ref{lem_partition}. We denote by 
$P_r:Y^\bdiscrete_\A \to Y^\bdiscrete_\A$ the projection defined by
\[
(P_r (\lambda_j)_{j \in I})_i = \left\{\begin{array}{ll} \lambda_i, & \mbox{if } i \in I_r,\\
0, &\mbox{otherwise.}\end{array}\right.
\]
Furthermore, define $S_r:Y^\bdiscrete_\A \to \Co Y_\A$, $r=1,\hdots,s$, 
to be the operator
$S_r (\lambda_i)_{i\in I} = f$ such that $V_g f(x_i) = \lambda_i$ for all $i \in I_r$.
The existence and boundedness of this operator follows from 
Theorem~\ref{thm_interpol}. Now we may define the following bounded linear
operators between $Y^\bdiscrete_\A$ and $\bigoplus_{r=1}^s \Co Y_\A$:
\[
\widehat{A}: Y^\bdiscrete_\A \to \bigoplus_{r=1}^s \Co Y_\A,\quad
\Lambda \mapsto (S_1 \circ P_1 \Lambda, \hdots, S_s \circ P_s \Lambda)
\]
 and
\[
\widehat{B}: \bigoplus_{r=1}^s \Co Y_\A \to Y^\bdiscrete_\A,\quad
(f_1,\hdots,f_s) \mapsto \sum_{r=1}^s P_r(V_g f(x_i)).
\]
It immediately follows that $\widehat{B} \circ \widehat{A} = \Id_{Y^{\bdiscrete}_\A}$.
Thus, $Y^\bdiscrete_\A$ is a retract of $\bigoplus_{r=1}^s \Co Y_\A$.
 \end{Proof}

This theorem allows to transform many questions about the invariant coorbit
spaces to questions about the corresponding sequence spaces. 
In order to formulate a particular result 
we recall that a class $\mathcal{J}$ of operators between arbitrary Banach spaces
is called an {\it operator ideal} if the following conditions are satisfied for the components
${\mathcal J}(E,F) := \mathcal{J} \cap \B(E,F)$, where $\B(E,F)$ denotes the space of 
bounded operators between Banach spaces $E,F$ 
(see also \cite[p.45]{Pietsch}).
\begin{itemize}
\item The identity operator $\Id_{\C}$ belongs to $\mathcal{J}$, where $\C$ is
identified with the one-dimensional Banach space.
\item If $S_1,S_2 \in {\mathcal J}(E,F)$ then $S_1 + S_2 \in \cJ(E,F)$.
\item If $T \in \B(E_0,E)$, $S \in \cJ(E,F)$ and $R \in \B(F,F_0)$ then 
$RST \in \cJ(E_0,F_0)$.
\end{itemize}
For instance, the compact operators form an operator ideal.
\begin{Theorem}\label{thm_coorbit_embed}
\begin{itemize}
\item[(a)] The inclusion mapping $J_\A :\Co Y_\A \hookrightarrow \Co Z_\A$  
between
invariant coorbit spaces is automatically continuous. The same holds for the
inclusion mapping 
$J^\bdiscrete_\A: Y^\bdiscrete_\A \hookrightarrow Z^\bdiscrete_\A$.
\item[(b)] $J_\A$ is compact if and only if $J^\bdiscrete_\A$ is compact.
\item[(c)] Let $\mathcal{J}$ be an operator ideal. Then $J_\A \in \mathcal{J}$
if and only if $J^\bdiscrete_\A \in \mathcal{J}$.
\end{itemize}
\end{Theorem}
\begin{Proof} The proof is exactly the same as the one of Theorem 9.4 in \cite{FG3}
and hence omitted.
\end{Proof}
We remark that this theorem may be used to estimate the entropy numbers
and approximation numbers of (compact) embedding operators of coorbit spaces.
For instance, the compact operators, whose entropy numbers (approximation numbers) 
are contained in some $\ell^p(\N)$-space form an operator ideal \cite{Pietsch}. 
So according to the above 
theorem it suffices to compute entropy or approximation numbers of
embeddings between sequence spaces, which is much easier than for function spaces, see
also the next section.

Let us now investigate compactness of embeddings for the important special case
$Y = L^p_m$. According to Lemma 4.3.1 in \cite{Rau_Diss} (see also \cite{Rau2})
the corresponding sequence space $Y^\bdiscrete_\A$ coincides with  $\ell^p_{m_p}$,
where $m_p(i) = m(x_i) |\A(x_iU)|^{1/p}$. The following lemmas give criterions
on the compactness of embeddings between certain weighted $\ell^p$-spaces.

\begin{Lemma}\label{lem_comp_embed} Let $1 \leq p\leq q\leq \infty$ and
$v,w$ be some positive weight functions on some infinite index set 
$I$. 
\begin{itemize} 
\item[(a)] $\ell^p_w(I)$ is continuously embedded into
$\ell^q_v(I)$ if and only if 
$\sup_{i \in I} v(i)/w(i) < \infty$. 
\item[(b)] The embedding $\ell^p_w(I) \hookrightarrow \ell^q_v(I)$ is compact 
if and only if the sequence $(v(i)/w(i))_{i \in I}$ is contained in $c_0(I)$, the space
of all sequences vanishing at $\infty$.
\end{itemize}
\end{Lemma}
\begin{Proof} The proof is straightforward and thus omitted. 
\end{Proof}

\begin{Lemma}\label{lem_pgq} 
Let $\infty \geq p > q \geq 1$ and $v,w$ be some positive 
weight functions on $I$.
Set $\beta := 1/q-1/p > 0$. If
$\sum_{i\in I} \left( \frac{v(i)}{w(i)}\right)^{1/\beta} < \infty$
then $\ell^p_w(I)$ is compactly embedded into $\ell^q_v(I)$.
\end{Lemma}
\begin{Proof} The assertion follows from a simple application of H\"older's inequality.
\end{Proof} 

Let us apply these lemmas to embeddings of coorbit spaces.

\begin{Theorem}\label{thm_compact} Suppose $\G$ is not compact. Let $v,m$ be two 
moderate invariant weight functions, and suppose $1 \leq p \leq q \leq \infty$.
Further, let $X=(x_i)_{i\in I}$ be some well-spread set with respect to $\A$ and 
$U=U^{-1}=\A(U)$ be some relatively compact neighborhood of $e \in \G$.
We define the sequence 
\begin{equation}\label{def_seqh}
h_i\,:=\,\frac{v_p(i)}{m_p(i)} \,=\,
\frac{v(x_i)}{m(x_i)|\A(x_i)U|^{1/p-1/q}},\qquad {i\in I}.
\end{equation}
\begin{itemize}
\item[(a)] $(\Co L^p_m)_\A$ is continuously embedded into $(\Co L^q_v)_\A$
if and only if $h \in \ell^\infty(I)$.
\item[(b)] The embedding in (a) is compact if and only if $h \in c_0(I)$.
\end{itemize}
\end{Theorem}
\begin{Proof} The sequence spaces associated to $L^p_m$ and $L^q_v$ are given by
$(L^p_m)^\bdiscrete_\A = \ell^p_{m_p}$ and $(L^q_v)_\A = \ell^q_{v_q}$ with 
(see \cite{Rau2} or \cite[Lemma 4.3.1]{Rau_Diss})
\begin{equation}\label{mp_weight}
m_p(i) \,=\, m(x_i) |\A(x_i U)|^{1/p} \quad\mbox{and}\quad v_q(i) \,=\, v(x_i) |\A(x_i U)|^{1/q}.
\end{equation}
By Lemma \ref{lem_comp_embed} $\ell^p_{m_p}$ is embedded into $\ell^q_{v_q}$ if and only if
the sequence $h$
is contained in $\ell^\infty$. Moreover, this embedding is compact if and only if $h \in c_0(I)$.
Hence, by Theorem~\ref{thm_seq_coorbit} we have $(\Co L^p_m)_\A \hookrightarrow (\Co L^q_v)_\A$ if and only
if $h \in \ell^\infty$. Moreover, this embedding is compact if and only if $h \in c_0(I)$.
\end{Proof}

It is also possible to get rid of the well-spread set $X$ in the previous theorem
as is shown by the next Lemma.

\begin{Lemma} With the same notation as in Theorem \ref{thm_compact} we define the
function
\begin{equation}\label{C0_cond}
H: \G \to \R,\quad x \,\mapsto\, \frac{v(x)}{m(x) |\A(x U)|^{1/p-1/q}}.
\end{equation}
Then we have 
$h \in \ell^\infty(I)$ if and only if $H \in L^\infty(\G)$ and $h \in c_0(I)$
if and only if $H \in C_0(\G)$.
\end{Lemma}
\begin{Proof}
Clearly, if the function $H$ in (\ref{C0_cond}) is contained in $L^\infty$, resp. $C_0(\G)$ 
then $h \in \ell^\infty$ resp. $h \in c_0$. This shows the "if"-part of (a) and (b).
The converse part follows from the moderateness of the functions
$v$ and $m$.
\end{Proof}

This theorem gives an abstract explanation to 
the phenomenon that restricting to functions that possess symmetry
may enforce compactness of embeddings. This fact was observed for the Besov and
Triebel-Lizorkin spaces of radial distributions recently in \cite{KLSS,SickS,LS,ST}.

Indeed, in typical situations the function
$x \mapsto |\A(xU)|$ has some growth. However, for the trivial subgroup $\A= \{e\}$
this function is constant by definition of the Haar measure. So by restricting 
to elements in $\Co Y$ which
are invariant under a suitable $\A$, we might be able to 
enforce the function in (\ref{C0_cond}) to belong to $C_0$, although the function
$v(x)/w(x)$ does not. 
We formulate this observation for a simple case in 
the next corollary.

\begin{Cor} Let $\G$ be non-compact. Suppose $\A$ is 
such that $x \mapsto |\A(x U)|^{-1}$ belongs to $C_0(\G)$ and let
$1 \leq p < q \leq \infty$. Further, let $m$ be some moderate invariant weight function. Then
$\Co L^p_m$ is embedded in $\Co L^q_m$, but not compactly. Restricting to $\A$-invariant
elements enforces compactness, i.e., the embedding 
$(\Co L^p_m)_\A \hookrightarrow (\Co L^q_m)_\A$ is compact.
\end{Cor}

Let us now consider the case $p>q$. 

\begin{Theorem}\label{thm_pgq} 
Let $\infty \geq p > q \geq 1$ and $v,m$ be two moderate invariant 
weight functions. Set $\beta:= 1/q - 1/p >0$.
If $\int_\G \left(\frac{v(x)}{m(x)}\right)^{1/\beta} dx < \infty$ then
$(\Co L^p_m)_\A$ is compactly embedded into $(\Co L^q_v)_\A$.
\end{Theorem}
\begin{Proof} The sequence spaces associated to $L^p_m$ and $L^q_v$ are given by
$(L^p_m)^\bdiscrete_\A = \ell^p_{m_p}$ and $(L^q_v)_\A = \ell^q_{v_q}$ with with $m_p$ and $v_q$ as in (\ref{mp_weight}). According
to Lemma \ref{lem_pgq} and Theorem \ref{thm_coorbit_embed} the assertion
follows from the following computation:
\begin{align}
\sum_{i\in I} \left(\frac{v_q(i)}{m_p(i)}\right)^{1/\beta}
\,&=\, \sum_{i \in I} \left(\frac{v(x_i)}{m(x_i)}|\A(x_i U)|^\beta\right)^{1/\beta} \,\asymp\, \sum_{i\in I} \left(\frac{v(x_i)}{m(x_i)}\right)^{1/\beta} \int_{\A(x_i U)} dx\notag\\
&\asymp\, \sum_{i\in I} \int_{\A(x_i U)} \left(\frac{v(x)}{m(x)}\right)^{1/\beta} dx
\,\asymp\, \int_\G \left(\frac{v(x)}{m(x)}\right)^{1/\beta} dx < \infty. \notag
\end{align}
Here the moderateness of $v,m$ and the finite overlap property
of the well-spread set $(x_i)_{i\in I}$ was used.
\end{Proof}

Since the condition in the previous theorem 
is independent of the automorphism group $\A$, restricting
to invariant elements does not give stronger results for embeddings -- in contrast
to the case $p<q$.

\subsection{Linear and Nonlinear approximation}

Let us now consider linear approximations of the form 
$\sum_{i\in N} \lambda_i \TPi(x_i) g$ of $f = \sum_{i\in I} \lambda_i \TPi(x_i) g$, 
where $N\subset I$ is a finite set. 
For some subspace $V$ of some Banach space $B$ we introduce the error
of approximation by
\begin{equation}\label{def_linapprox_error}
e(f,V,B) \,:=\, \inf_{g \in V} \|f-g|B\| \qquad \mbox{for } f \in B.
\end{equation}
Furthermore, for some bounded sequence $\lambda=(\lambda_i)_{i\in I}$ we
denote by $(s_n(\lambda))_{n\in \N}$ its non-increasing rearrangement, i.e.,
\begin{align}\label{def_incre}
s_n(\lambda) \,:=\, \inf \left\{\sigma \geq 0, \#\{i: |\lambda_i|\geq \sigma\} < n \right\}.
\end{align}
Clearly, it holds $s_n(\lambda) \geq s_{n+1}(\lambda)$ for all $n \in \N$. 

\begin{Theorem}\label{thm_linapprox} Let $1\leq p \leq q \leq \infty$ and $v,m$ 
be some moderate weight functions on $\G$ such that
the condition in Theorem (b) is satisfied, i.e., the embedding 
$(\Co L^p_m)_\A \hookrightarrow (\Co L^q_v)_\A$
is compact.
Let $X=(x_i)_{i \in I}$ and $g \in \BwA$ 
such that $\{\TPi(x_i) g\}_{i\in I}$ is an atomic 
decomposition of
$(\Co L^p_m)_\A$ and $(\Co L^q_v)_\A$ as in Theorem \ref{thm_compact}. 
We denote by $(h_i)_{i\in I}$ the sequence defined in (\ref{def_seqh}).
Let $\tau : \N \to I$ realize the non-increasing rearrangement of $h$, i.e.,
$s_n(h) = h_{\tau(n)}$, and let $V_n = \spann \{\TPi(x_{\tau(j)})g, j=1,\hdots,n\}$.
Then for all $n \in \N$
\[
e(f,V_n,(\Co L^q_v)_\A) \,\leq\, C s_{n+1}(h) \|f|(\Co L^p_w)_\A\| \quad\mbox{for all }
f \in (\Co L^p_m)_\A.
\]
\end{Theorem}
\begin{Proof}  Let $\lambda_i = \langle f, e_i\rangle$ with $e_i \in (\H^1_w)_\A$ as in
Theorem \ref{thm_ab}. It holds $(L^p_m)^\ddiscrete = \ell^p_{m^{(p)}}$ and $(L^q_v)^\bdiscrete = \ell^q_{v^{(q)}}$ with $m^{(p)}(i) = m(x_i) |\A(x_i U)|^{1/p-1}$ and 
$v^{(q)} =  v(x_i) |\A(x_i U)|^{1/q-1}$ by Lemma 4.3.1 in \cite{Rau_Diss}.
Since $\|(\lambda_i)_{i\in I}|\ell^q_{m^{(q)}}\|$ is an 
equivalent norm on $(\Co L^p_m)_\A$ by Theorem \ref{thm_ab}
we obtain for all $f \in (\Co L^p_m)_\A$
\begin{align}
& \|f - \sum_{j=1}^n \lambda_{\tau(j)} \TPi(x_{\tau(j)})g|(\Co L^q_v)_\A\|
\,=\, \|\sum_{i\in I} \lambda_i \TPi(x_i)g - 
\sum_{j=1}^n \lambda_{\tau(j)} \TPi(x_{\tau(j)}) g|(\Co L^q_v)_\A\|\notag\\
& =\, \|\sum_{j=n+1}^\infty \lambda_{\tau(j)} \TPi(x_{\tau(j)}) g|(\Co L^q_v)_\A\|
\leq C \left(\sum_{j=n+1}^\infty (|\lambda_{\tau(j)}| v^{(q)}(\tau(j)))^q\right)^{1/q}\notag\\
&\leq\, C \sup_{j>n} \frac{v(x_{\tau(j)}) |\A(x_{\tau(j)}
U)|^{1/q-1}}{w(x_{\tau(j)})|\A(x_{\tau(j)}U)|^{1/p-1}} \left(\sum_{j=n+1}^\infty
(|\lambda_{\tau(j)}|^p m^{(p)}(\tau(j))^p\right)^{1/p}\notag\\
&\leq\, C' h_{\tau(n+1)} \|f|(\Co L^p_m)_\A\|.\notag
\end{align}
Hereby, we used also Lemma \ref{lem_comp_embed}. This yields the claim.
\end{Proof} 

\begin{Remark}\label{rem_lin} 
This theorem shows that using elements $\TPi(x) g$ instead of the elements
$\pi(y) g$ for approximating $f \in \Co (L^p_w)_\A$ gives an advantage. Indeed, in 
typical situations
$|\A(xU)|$ is a growing function so that the sequence $h = h_\A$ in (\ref{def_seqh}) 
decreases faster than $h_{\{e\}}$, the one for the trivial automorphism group $\{e\}$. 
This means that the error of linear approximation (measured in the $\Co L^q_v$-norm) 
with elements $\TPi(x)g$ decreases faster
than the one of approximation with elements $\pi(y) g$. Moreover, if $v=w$ and the function 
$x\mapsto |\A(x U)|^{-1}$ is contained in $C_0(\G)$ then the statement in Theorem 
\ref{thm_linapprox}
is a significant improvement for invariant $f$'s since for the approximation
with $\pi(y)g$ we only know that the error (measured again in $\Co L^q_v$) 
converges to $0$ (with no information about the speed
of convergence and provided $q<\infty$) while we have a more concrete error estimate for the 
approximation with elements $\TPi(x) g$. 
\end{Remark} 

Let us finally discuss non-linear approximation. 
Let $(x_i)_{i\in I}$ be some well-spread set with respect to $\A$ and $g \in \BwA$ 
such that $\{\TPi(x_i)g\}_{i\in I}$ forms an atomic decomposition.
We denote by
\[
\sigma_n(f,\Co Y_\A) \,:=\, \inf_{N \subset I, \#N \leq n} \|f - \sum_{i \in N} \lambda_i \TPi(x_i) g|\Co Y_\A\|
\]
the error of best $n$-term approximation. Hereby, the infimum is also taken over all possible
choices of coefficients $\lambda_i$. 
Our task is to find a class of elements for which this error has a certain decay
when $n$ tends to $\infty$.

The following lemma, which is taken from \cite{GS}, is useful for this task.

\begin{Lemma}\label{lem_nonlinear} 
Let $b=(b_k)_{k\in \N}$ be some non-increasing sequence of positive
numbers. Set $\sigma_{n,q}(b) = (\sum_{k=n}^\infty b_k^q)^{1/q}$ and for $p,q > 0$ set $\alpha = 1/p-1/q$.
Then for $0<p<q\leq \infty$ we have
\[
2^{-1/p} \|b|\ell^p(\N)\| \,\leq\, \left(\sum_{n=1}^\infty (n^\alpha \sigma_{n,q}(b))^p
\frac{1}{n}\right)^{1/p} \,\leq\, C \|b|\ell^q(\N)\|.
\] 
\end{Lemma}

\begin{Theorem}\label{thm_nonlinapprox} 
Let $m,v$ be some $w$-moderate weight functions on $\G$, let $1\leq p < q \leq \infty$
and define $\alpha = 1/p-1/q$.
Let $(x_i)_{i\in I}$ be some well-spread set with respect to $\A$ 
such that $\{\TPi(x_i)g\}_{i\in I}$ is an atomic decomposition
of $(\Co L^p_m)_\A$. Assume further that the function $H$ defined in (\ref{C0_cond})
is contained in $L^\infty$, or equivalently $(\Co L^p_m)_\A \hookrightarrow (\Co L^q_v)_\A$.
Then 
\begin{equation}\label{best_nterm_decay}
\left(\sum_{n=1}^\infty \frac{1}{n}(n^\alpha \sigma_n(f,(\Co L^q_v)_\A))^p\right)^{1/p}
\,\leq\, C \|f|(\Co L^p_m)_\A\| \quad \mbox{for all } f\in (\Co L^p_m)_\A.
\end{equation}
\end{Theorem}
\begin{Proof} Let $f = \sum_{i \in I} \lambda_i \TPi(x_i)g$ with 
$\lambda_i = \langle f,e_i\rangle$ be an expansion of 
$f \in (\Co L^p_m)_\A$ in terms of the atomic decomposition.  
Further let $v^{(q)}(i)=v(x_i) |\A(x_iU)|^{1/q-1}$ and $b(i) := |\lambda_i| v^{(q)}(i)$. 
Let $\tau : \N \to I$ be a
bijection such that
$b(\tau(s)) \geq b(\tau(s+1))$ for all $s \in \N$.  
We obtain
\begin{align}
\sigma_n(f, (\Co L^q_v)_\A) \,&\leq\, 
\|\sum_{s=n+1}^\infty \lambda_{\tau(s)} \TPi(x_{\tau(s)})g|(\Co L^q_v)_\A\|
\,\leq\, C (\sum_{s=n+1}^\infty |\lambda_{\tau(s)} v^{(q)}(\tau(s))|^q)^{1/q}\notag\\
&=\, C \sigma_{n,q}(b).
\end{align}
By Lemma \ref{lem_nonlinear} and Theorem \ref{thm_ab} we deduce
\begin{align}
&\left(\sum_{n=1}^\infty \frac{1}{n}(n^\alpha \sigma_n(f,(\Co L^q_v)_\A))^p\right)^{1/p}
\,\leq\, C 
\left(\sum_{n=1}^\infty \frac{1}{n}(n^\alpha 
\sigma_{n,q}(b)^p\right)^{1/p}\notag\\
&\leq\, C' \left(\sum_{i \in I} (|\lambda_i| v^{(q)}(i)|)^p\right)^{1/p}
\,=\, C' \left(\sum_{i\in I} |\lambda_i|^p v(x_i)^p |\A(x_i U)|^{(1/q-1)p}\right)^{1/p}\notag\\
&=\, C' \left(\sum_{i\in I} |\lambda_i|^p (v(x_i) |\A(x_i U)|^{1/q-1/p})^p |\A(x_i U)|^{p(1/p-1)}\right)^{1/p}\notag\\
&\leq \, C'' \left(\sum_{i \in I} |\lambda_i|^p (m(x_i) |\A(x_i U)|^{1/p-1})^p\right)^{1/p}
\,\leq\, C''' \|f|(\Co L^p_m)_\A\|.\notag
\end{align}
\end{Proof}
Note that (\ref{best_nterm_decay}) implies 
\[
\sigma_n(f, (\Co L^q_v)_\A) \leq C n^{-\alpha} \quad  \mbox{for all } f \in (\Co L^p_m)_\A.
\]
We remark that once again this theorem shows that using the atomic decomposition 
$\{\TPi(x_i)g\}_{i\in I}$ instead
of the atomic decomposition $\{\pi(y_j)g\}_{j \in J}$ 
for approximating elements in $(\Co L^p_m)_\A$ is advantageous.
Indeed, the class of functions $f$ for which the error of best $n$-term approximation
with $\{\TPi(x_i)g\}_{i\in I}$
has some prescribed decay (i.e. for which an estimate as in (\ref{best_nterm_decay}) holds)
is stricly larger than the one for $\{\pi(y_j)g\}_{j\in J}$ in typical situations.
Indeed, if $x \mapsto |A(xU)|$ is an unbounded function 
then we may have $(\Co L^p_m)_\A \subset (\Co L^q_v)_\A$ although 
$\Co L^p_m \not\subset \Co L^q_v$. 

Let us note further that Theorem \ref{thm_nonlinapprox} does not give a characterization of all 
elements $f$ satisfying the dacay condition (\ref{best_nterm_decay}), i.e., we cannot provide a converse
inequality of (\ref{best_nterm_decay}) (a Bernstein inequality).
This is due to the fact
that we do not work with a basis but with a frame. 
In fact, it is a difficult (and open) problem 
to find classes of frames for which Bernstein inequalities hold.

\section{Radial Time Frequency Analysis}
\subsection{Short Time Fourier Transform and Modulation Spaces}

We will now apply the abstract results from the previous section to time frequency analysis
of radial functions.

Let $\HH_d := \R^d \times \R^d \times \TT$ denote the (reduced) Heisenberg group with group law
\[
(x,\omega,\tau) (x',\omega',\tau') ~=~ 
(x+x',\omega + \omega', \tau \tau' e^{\pi i(x'\cdot \omega - x\cdot \omega')}).
\]
The Heisenberg group is unimodular and has Haar measure 
\[
\int_{\HH_d} f(h) dh = \int_{\R^d} \int_{\R^d} \int_0^1 f(x,\omega,e^{2\pi i t}) dt d\omega dx.
\]
The Schr\"odinger representation $\rho$ acting on $\H=L^2(\R^d)$ is described as follows.
Let 
\[
T_x f(t) \,:=\, f(t-x), \quad \mbox{and}\quad M_\omega f(t) \,=\, e^{2\pi i \omega \cdot t} f(t), \quad
x,\omega, t \in \R^d,
\]
denote the translation and modulation operator on $L^2(\R^d)$. Then $\rho$ is defined by
\[
\rho(x,\omega,\tau) \,:=\, \tau e^{\pi i x\cdot \omega} T_x M_\omega
\,=\, \tau e^{-\pi i x \cdot \omega} M_\omega T_x.
\]
It is well-known that this is an irreducible unitary and square-integrable representation of $\HH_d$.
The corresponding voice transform is essentially the short time Fourier 
transform:
\begin{align}
V_g f(x,\omega,\tau) \,=&\, \langle f, \rho(x,\omega,\tau) g\rangle_{L^2(\R^d)} 
\,=\, \overline{\tau} 
\int_{\R^d} f(t) \overline{e^{-\pi i x\cdot \omega} M_\omega T_x g(t)} dt\notag\\
\label{Schroed_STFT}
=&\, \overline{\tau} e^{\pi i x\cdot \omega} 
\int_{\R^d} f(t) \overline{g(t-x)} e^{-2\pi i t \cdot \omega} dt 
\,=\, \overline{\tau} e^{\pi i x\cdot \omega} \STFT_g f(x,\omega).
\end{align}

The automorphisms of $\R^d \times \R^d$ that extend to automorphisms
of $\HH_d$ are given by the elements of the symplectic group $Sp(d)$.
The latter is defined as the subgroup of $GL(2d,\R)$ leaving invariant 
the symplectic form
$[(x,\omega),(x',\omega')]:= x'\cdot \omega - x\cdot \omega'$ (see \cite{GrBook}).
A compact subgroup of $Sp(d)$ is given by
\begin{equation}\label{auto_SOd}
\A \,:=\, \left\{ \left(\begin{array}{cc} A & 0 \\ 0 & A \end{array}\right)
~|~ A \in SO(d) \right\} \cong SO(d).
\end{equation}
An element 
$A \in SO(d) \cong \A$ acts on $\HH_d$ by 
$A (x,\omega,\tau) = (Ax, A\omega, \tau)$. 

We choose the natural 
representation $\sigma$ of $SO(d)$ on $L^2(\R^d)$ given by
$\sigma(A) f(t)=f(A^{-1} t)$ for $A \in SO(d), t \in \R^d$.
Using the orthogonality of $A \in SO(d)$ we obtain
\begin{align}
\rho(Ax,A\omega,\tau) \sigma(A) f(t)
\,=&\; \tau e^{-\pi i (Ax\cdot A\omega)} e^{2\pi i A\omega \cdot t}
f(A^{-1}(t-A x))\notag\\
\,=&\; \tau e^{-\pi i (x\cdot \omega)} e^{2\pi i \omega \cdot A^{-1} t}
f(A^{-1}t - x)
\,=\, \sigma(A) \rho(x,\omega,\tau)f(t).\notag
\end{align}
Thus, condition (\ref{Cond_pi_sigma}) is satisfied. Clearly, it holds
$\H_\A = L^2_{rad}(\R^d)$, the space of radially symmetric $L^2$-functions.
In the sequel we assume $d\geq 2$.

The operator $\tilde{\rho}$ as defined in (\ref{def_TPi}) reads
\begin{align}
\tilde{\rho}(x,\omega,\tau) f(t) \,&=\, \int_{SO(d)} \rho(Ax,A\omega,\tau)f(t) dA 
\,=\, \tau e^{\pi i x \cdot \omega} \int_{SO(d)} e^{2\pi i A \omega\cdot t} f(t-Ax) dA\notag\\
&=:\, \tau e^{\pi i x \cdot \omega} \Omega(x,\omega)f(t),\quad (x,\omega,\tau) \in \HH_d.\notag
\end{align}

\begin{figure}
\parbox[t]{6.8cm}{\includegraphics[width=6.8cm]{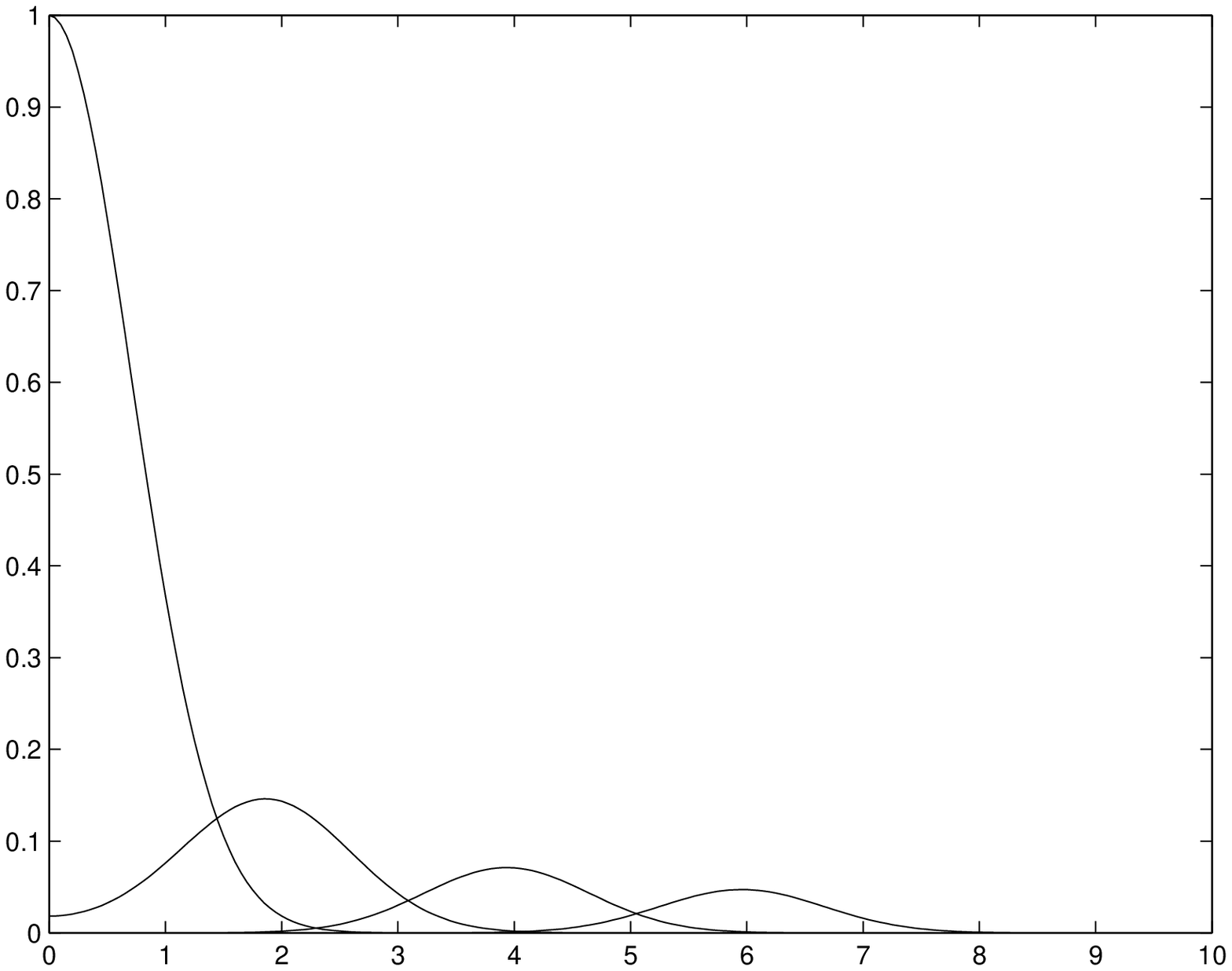} 
\center{$\Omega(r,0,1)f$ for $r=0,2,4,6$}}
\parbox[t]{6.8cm}{\includegraphics[width=6.8cm]{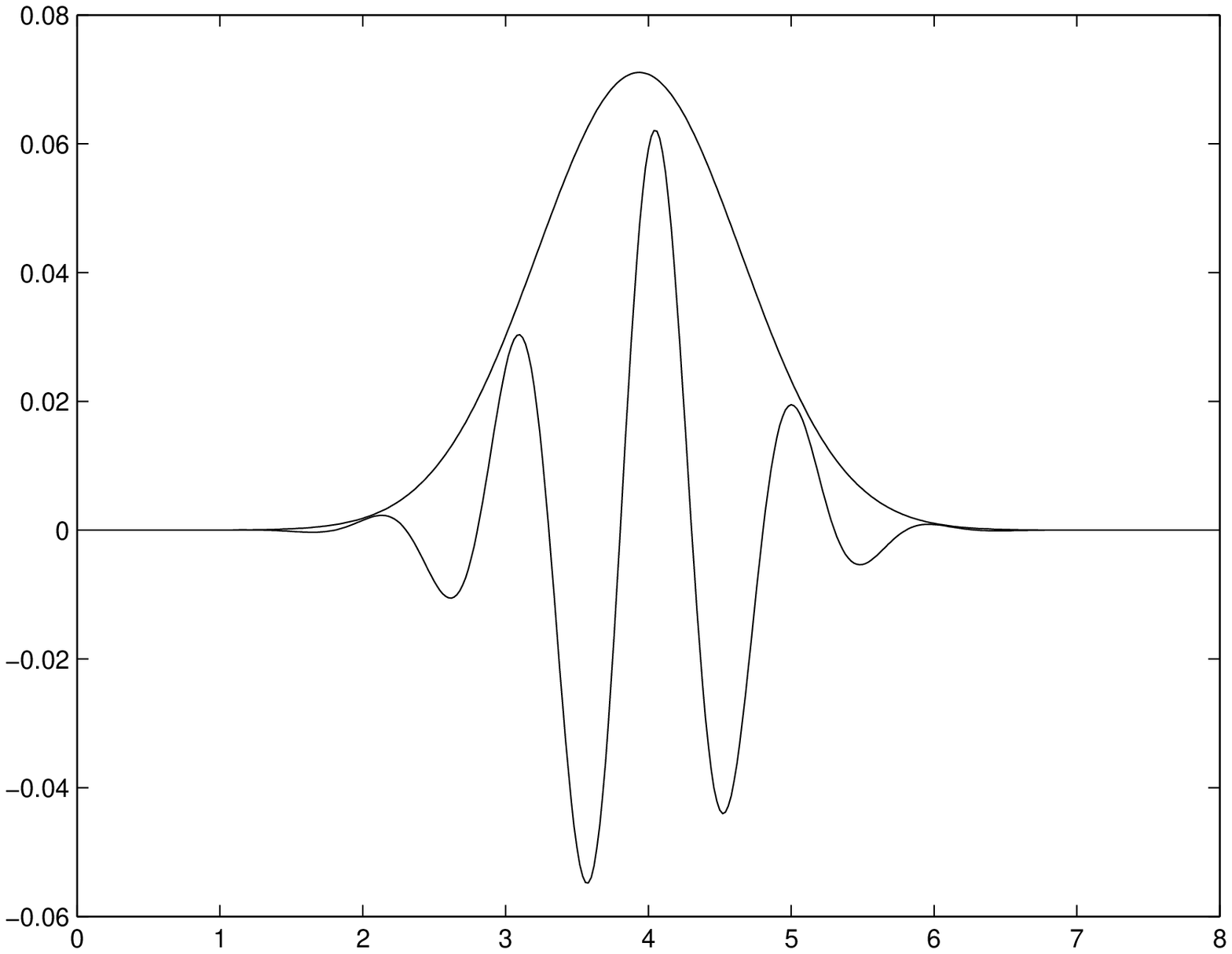} 
\center{$\Re(\Omega(4,s,1)f)$ for $s=0,1$}}
\caption{$\Omega(r,s,\cos\alpha)f$ for a Gaussian $f(x) = e^{-|x|^2}$ for
different values of $r,s,\alpha$}\label{plots1}
\end{figure}

It has been shown in \cite{Rau1,Rau_Diss} that $\Omega$ may be expressed as
\begin{align}
&\Omega(x,\omega)f(t) \,=\,\Omega(r,s,\cos \alpha) f_0(\theta) \notag\\
&=\, \frac{|S^{d-2}|}{|S^{d-1}|} 
\int_0^\pi f_0(\sqrt{\theta^2 - 
2 r\theta \cos \phi+ r^2}) e^{2\pi i \theta s \cos \alpha \cos \phi}
\B_{d-1}(\theta s \sin \alpha \sin \phi) \sin^{d-2} \phi\, d\phi\notag
\end{align}
where $r = |x|, s = |\omega|, 
x \cdot\omega = rs \cos \alpha, \theta = |t|$ and $f_0 : [0,\infty) \to \C$ is 
such that $f(x) = f_0(|x|)$.
Here, $\B_{d-1}$
denotes the spherical Bessel function defined by 
\begin{equation}\label{sph_Bessel}
\B_d(t) ~=~ \frac{1}{|S^{d-1}|} \int_{S^{d-1}} 
e^{2\pi i t \eta \cdot \xi} dS(\xi),\quad \eta \in S^{d-1}.
\end{equation}
(independent of the choice of $\eta \in S^{d-1}$). 

\begin{figure}
\parbox[t]{6.8cm}{\includegraphics[width=6.8cm]{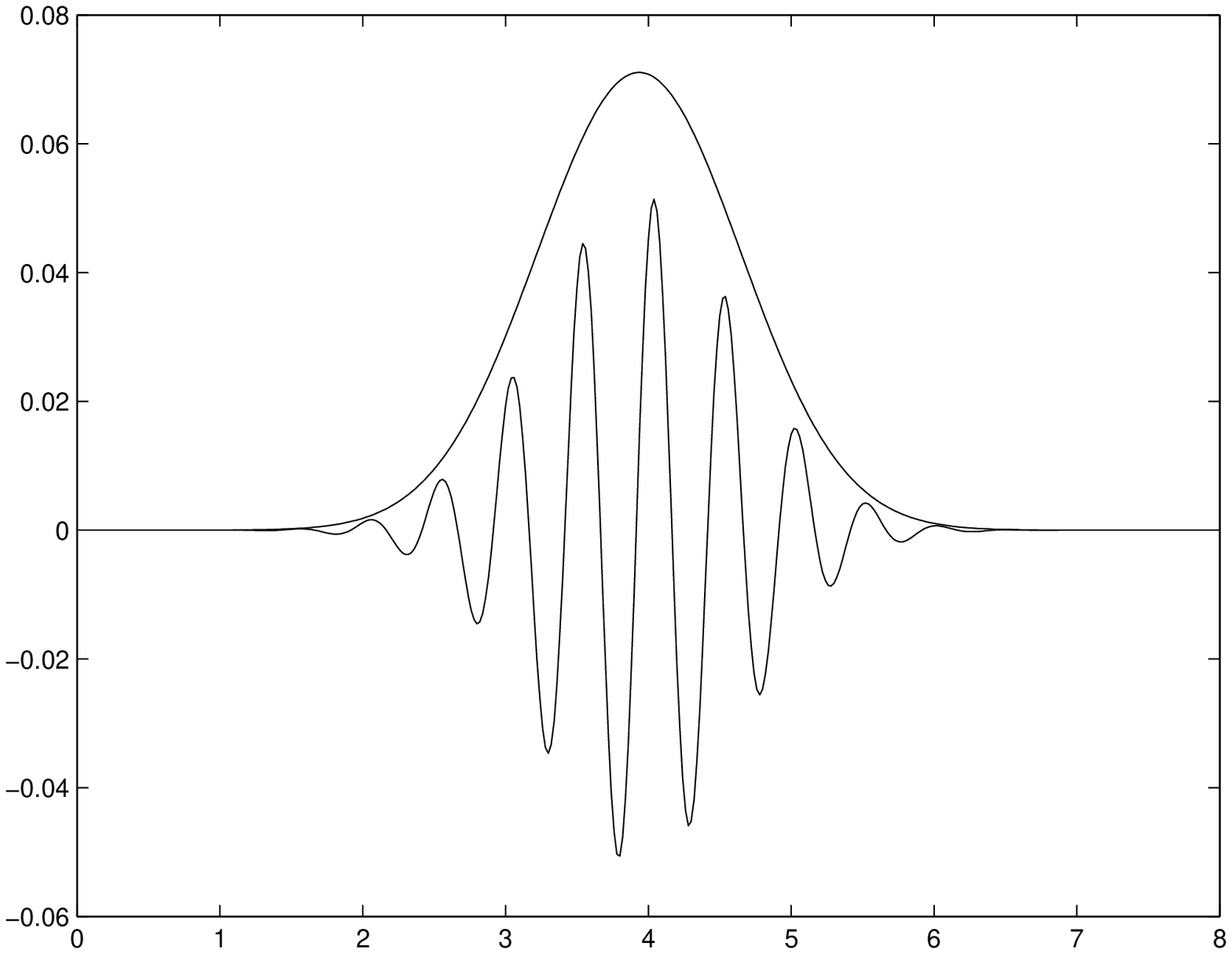} 
\center{$\Re(\Omega(4,s,1)f)$ for $s=0,2$}}
\parbox[t]{6.8cm}{\includegraphics[width=6.8cm]{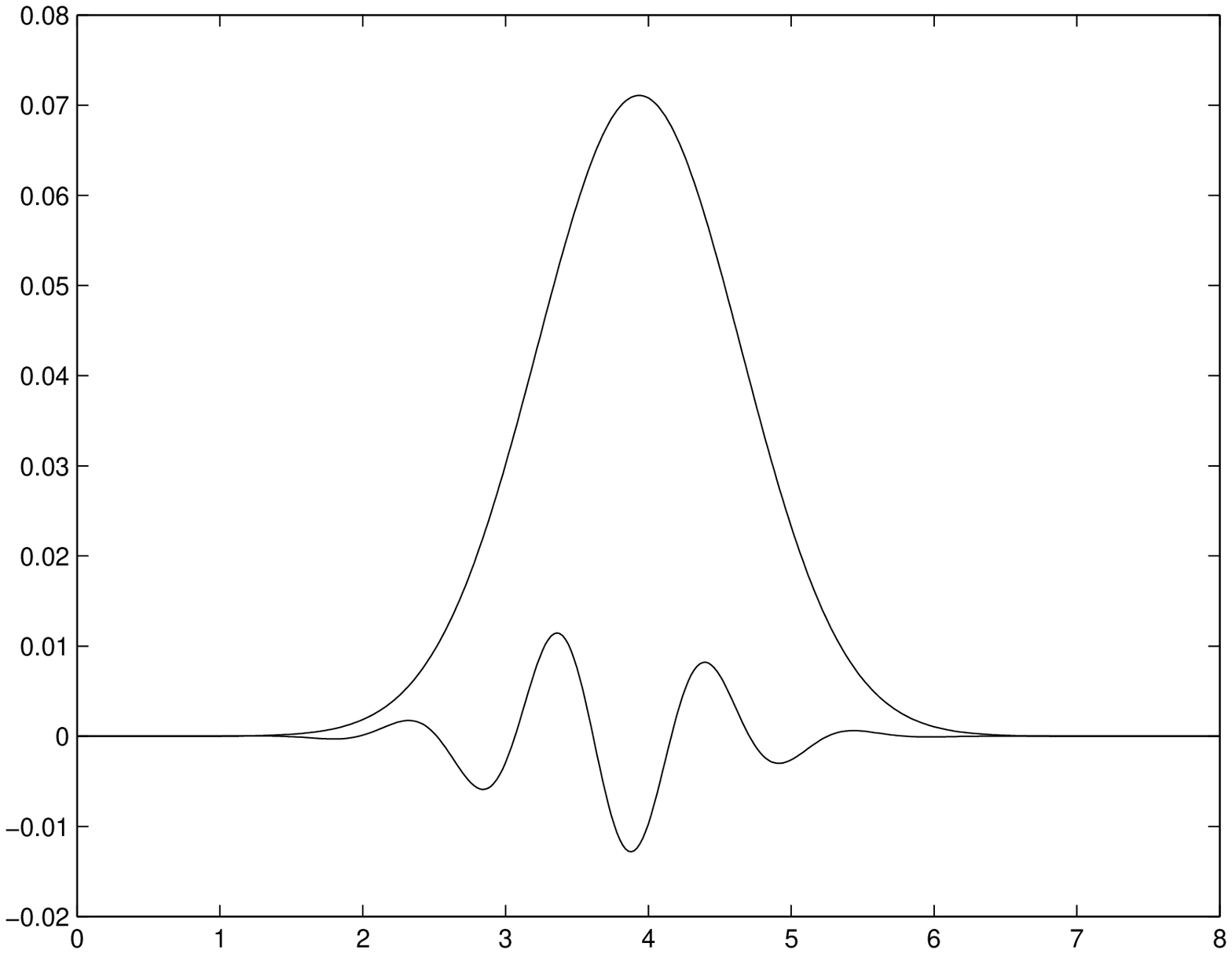} 
\center{$\Re(\Omega(4,s,\cos \alpha)f)$ for $(s,\alpha)=(0,0),(1,\pi/6)$}}
\caption{$\Omega(r,s,\cos\alpha)f$ for a Gaussian $f(x) = e^{-|x|^2}$ for
different values of $r,s,\alpha$}\label{plots2}
\end{figure}

Note that the spherical Bessel function (\ref{sph_Bessel}) can be expressed by means 
of the Bessel function $J_\alpha$ of
the first kind
\begin{equation*}
\B_d(t) \,=\,  \Gamma(\alpha+1) (\pi t)^{-\alpha} J_\alpha(2\pi t),\quad \alpha = \frac{d-2}{2}.
\end{equation*}
In particular, we have
\begin{align}
\B_1(t) \,=\, \cos(2\pi t),
\qquad
\B_2(t) \,=\, J_0(2\pi t)
\qquad\mbox{and}\qquad
\B_3(t) \,=\, \frac{\sin(2\pi t)}{2\pi t}.\notag
\end{align}
Figures \ref{plots1}, \ref{plots2} illustrate $\Omega(r,s,\cos\alpha)f$ for a Gaussian $f$ and $d=2$
for certain values of $r,s,\alpha$. We have always plotted the real part
of the functions (as functions of $|x|$, $x \in \R^2$.)

For more details about the short time Fourier transform of radial functions we refer to
\cite{Rau1,Rau_Diss}. 

Let us now introduce the modulation spaces on $\R^d$. We consider nonnegative
continuous weight functions $m$ on $\R^d \times \R^d$ that satisfy
\[
m(x+y,\omega + \xi) \,\leq\, C (1 + |x|^2 + |\omega|^2)^{a/2} m(y,\xi), \quad (x,\omega),(y,\xi) \in \R^d \times \R^d.
\]
for some constants $C>0,a \geq 0$. 
This means that $m$ is a moderate function with respect to $(x,\omega)\mapsto (1+|x|^2+|\omega|^2)^{a/2}$,
see also \cite[Chapter 11.1]{GrBook}.
A typical choice is 
\begin{equation}\label{def_ms}
m_s(x,\omega) \,=\, (1+|\omega|)^{s}, \quad s \in \R.
\end{equation}
Now let $g$ be some non-zero Schwartz function on $\R^d$, for instance a Gaussian.
The short time Fourier transform $\STFT_g$ extends to the space $\cS'(\R^d)$ of tempered distributions
in a natural way. 
Given $1\leq p,q\leq \infty$ and $m$ as above the modulation space $M^{p,q}_m$ is defined
as the collection of all distributions $f \in \cS'(\R^d)$ for which the norm
\begin{align}\label{Def_Mod}
\|f|M^{p,q}_m\| \,=\, \left(\int_{\R^d}\left(\int_{\R^d} |\STFT_g f(x,\omega)|^p m(x,\omega)^p dx\right)^{q/p}d\omega\right)^{1/q} 
\end{align}
is finite. We will sometimes 
restrict to the case $p=q$ in the sequel and denote 
$M^p_m = M^{p,p}_m$.

Since $|V_gf(x,\omega,\tau)| = |\STFT_g f(x,\omega)|$ 
we can easily identify the modulation
spaces with coorbit spaces, i.e.
\[
M^{p,q}_m(\R^d) \,=\, \Co L^{p,q}_m(\HH_d),
\]
where $m$ is extended to $\HH_d$ in a trivial way 
by $m(x,\omega,\tau) = m(x,\omega)$ and $L^{p,q}_m$ is a weighted
mixed norm space whose definition is obvious from (\ref{Def_Mod}).

Moreover, the elements invariant under $\A$ are clearly the 
radial distributions, hence
\[
(\Co L^{p,q}_m)_{SO(d)} \,=\, (M^{p,q}_m)_{rad}(\R^d) \,=\, \{f \in M^p_m, f \mbox{ is radial}\}.
\]

\subsection{Radial Gabor Frames}

We are interested in atomic decompositions and Banach frames
of the form $\{\Omega(x_i,\omega_i) g\}_{i\in I}$ of $M^{p,q}_m$, i.e., 
we want to apply Theorem \ref{thm_ab} to this particular situation.
To this end we need a covering of $\R^d \times \R^d$ by sets of the form $\A((x_i,\omega_i) + U)$ with
$\A \cong SO(d)$ being the automorphism group in (\ref{auto_SOd}) 
and $U=U^{-1}=\A(U)$ 
some relatively compact neighborhood of $0 \in \R^d \times \R^d$. 

Given $a,b > 0$ let $U_{a,b}:= B(0,a) \times B(0,b) \subset \R^d \times \R^d$, where
$B(x,r)$ denotes the closed ball in $\R^d$ of radius $r$ centered at $x$. Clearly,
$U_{a,b}$ is invariant under $SO(d)$. For $(x,\omega) \in \R^d \times \R^d$ we consider
the set 
\begin{align}
V_{a,b}(x,\omega)\,:=&\,
SO(d) ((x,\omega) + U_{a,b})\notag\\ 
\,=&\, \{(A(x+y),A(\omega +\xi)),\, A \in SO(d), y \in B(0,a), \xi \in B(0,b) \}.\notag
\end{align}
By construction this set is again invariant under $SO(d)$. It follows from 
Lemma 4.1 in \cite{Rau1} that $V_{a,b}(x,\omega)$ depends only on
$|x|$, $|\omega|$ and $x\cdot \omega$.

In \cite{Rau_Diss} a relatively separated covering of $\R^d \times \R^d$ of the form
$\{V_{a,b}(x_i,\omega_i)\}_{i\in I}$ was constructed. Let
\begin{align}
N(j,k) \,:=&\, 
\left\lceil\frac{\pi}{4} \left(\arctan 
\frac{ k \left(3+\frac{3}{2j} - (\frac{3}{4j})^2\right)^{1/2} + j \left(3+\frac{3}{2k} - (\frac{3}{4k})^2\right)^{1/2}}
{(jk + \frac{1}{2}(j + k) - \frac{3}{8}(\frac{j}{k} + \frac{k}{j})+1)}\right)^{-1} \right\rceil,
\quad j,k \in \N,\notag\\
N(0,k) \,:=&\, N(j,0) \,:=\, 0, \qquad j,k \in \N_0.\notag
\end{align}
We note that asymptotically these numbers behave like
\begin{equation}\label{Njk_asymp}
N(j,k) \,\asymp\, \frac{\pi}{4\sqrt{3}}  \frac{jk}{j + k}.
\end{equation}
We further define 
\begin{align}
\theta_{j,k}^\ell \,:=&\, \cos \alpha_{j,k}^\ell \,:=\, \sin \frac{\pi \ell}{2 N(j,k)},
\quad j,k \in \N, \ell = -N(j,k),\hdots, N(j,k),\notag\\
\theta_{0,k}^0 \,:=&\, \theta_{j,0}^0 \,=\, 1, \quad j,k \in \N_0, \notag
\end{align}
and, finally, with some fixed unit vectors $\eta,\zeta$ with $ \eta \cdot \zeta = 0$ we let
\[
(x_{j,k,\ell},~\omega_{j,k,\ell}) \,:=\, (a j \eta,~
b k (\cos (\alpha_{j,k}^\ell) \eta + \sin (\alpha_{j,k}^\ell) \zeta)),\qquad
(j,k,\ell) \in I,
\]  
where the index set is given by 
\[
I  \,:=\, \{(j,k,\ell),\, j,k \in \N_0, \ell=-N(j,k),\hdots, N(j,k) \}.
\] 
We further denote by $X_{a,b}$ the collection of all points 
$\{(x_{j,k,\ell},\omega_{j,k,\ell}),\, (j,k,\ell) \in I\}$ emphasizing
the dependence on the parameters $a,b$.

\begin{Theorem}(\cite[Theorem 4.7.4]{Rau_Diss}) For all $a,b > 0$ 
the point set 
$X_{a,b}$ is a $U_{a,b}$-dense,
well-spread set with respect to $SO(d)$. In particular, it holds
\begin{equation}\label{cover_invariant}
\R^d \times \R^d \,=\, \bigcup_{(j,k,\ell) \in I} V_{a,b}(x_{j,k,\ell},\omega_{j,k,\ell}).
\end{equation}
\end{Theorem}

Furthermore, we need some estimation of the Lebesgue-measure of $V_{a,b}(x,\omega)$ in $\R^d \times \R^d$.
The following theorem follows from Lemma 4.7.3 in \cite{Rau_Diss}, see also \cite[p.145]{Rau_Diss}.

\begin{Theorem}\label{thm_measVab} Let 
\begin{align}
\mu_{j,k,\ell} \,:=& \, (j+k) \left(jk \cos \frac{\pi \ell}{2N(j,k)}\right)^{d-2},
(j,k,\ell) \in I,\, \ell \neq \pm N(j,k), j,k \geq 1,\notag\\
\mu_{j,k,\pm N(j,k)} \,:=&\, j^{d-1} + k^{d-1} + 1,
\quad j,k \in \N_0.\notag 
\end{align}
Then there are constants $C_1=C_1(a,b), C_2 = C_2(a,b)$ such that
\[
C_1 \mu_{j,k,\ell} \,\leq\, |V_{a,b}(x_{j,k,\ell},\omega_{j,k,\ell})| \,\leq\, C_2 \mu_{j,k,\ell}.
\]
\end{Theorem}

With this explicit covering of $\R^d \times \R^d$ at hand we
may apply Theorem \ref{thm_ab} to obtain radial Banach frames for coorbit spaces,
see also Theorem 4.7.6 in \cite{Rau_Diss}.

\begin{Theorem}\label{rad_mod_frame} Let $w$ be a submultiplicative
$SO(d)$-invariant weight 
and suppose $g \in (M^{1,1}_w)_{rad}$. Then there exist constants $a_0, b_0>0$ such that for all
positive $a<a_0, b<b_0$, $1\leq p,q \leq \infty$ and all $w$-moderate $SO(d)$-invariant 
weights $m$ the following holds:
\begin{itemize}
\item[(i)] The set of functions
\begin{align}
\{g_{j,k,\ell}\}_{(j,k,\ell)\in I} 
\,:=&\, \left\{\Omega\left(aj,bk,\sin \frac{\pi \ell}{2N(j,k)}\right) g,~ (j,k,\ell) \in I \right\}
\notag
\end{align}
forms a Banach frame of $(M^{p,q}_m)_{rad}$ with corresponding sequence
space $(L^{p,q}_m)_{SO(d)}^\bdiscrete$.
\item[(ii)] The functions $\{g_{j,k,\ell}\}_{(j,k,\ell)}$ form an atomic
decomposition of $(M^{p,q}_m)_{rad}$ with corresponding sequence
space $(L^{p,q}_m)_{SO(d)}^\ddiscrete$. In particular, any $f \in (M^{p,q}_m)_{rad}$ has an expansion
\[
f \,=\, \sum_{j = 0}^\infty \sum_{k=0}^\infty \sum_{\ell=-N(j,k)}^{N(j,k)} \lambda_{j,k,\ell}(f)
\,\Omega\left(aj,bk, \sin \frac{\pi \ell}{2N(j,k)}\right) g
\]
with norm-convergence whenever $1\leq p,q < \infty$, and weak-$*$ convergence otherwise. 
Moreover, we have the norm-equivalence
\[
\|(\lambda_{j,k,\ell}(f))_{(j,k,\ell) \in I}|(L^{p,q}_m)_{SO(d)}^\ddiscrete\| 
\,\asymp\, \|f|(M^{p,q}_m)_{rad}\|.
\]
\end{itemize}
\end{Theorem}

In the special case $p=q$ we have the following nice description of the sequence spaces. For some
moderate weight function $m$ let 
\[
m_{j,k,\ell} \,:=\, m(x_{j,k,\ell},\omega_{j,k,\ell}).
\]
Then the following norms are equivalent, 
\begin{align}
\|(\lambda_{j,k,\ell})|(L^{p}_m)_{SO(d)}^\bdiscrete\|
\,\asymp\, \left(\sum_{(j,k,\ell) \in I} |\lambda_{j,k,\ell}|^p m_{j,k,\ell}^p\, \mu_{j,k,\ell} \right)^{1/p},
\notag\\
\|(\lambda_{j,k,\ell})|(L^{p}_m)_{SO(d)}^\ddiscrete\|
\,\asymp\, \left(\sum_{(j,k,\ell) \in I} |\lambda_{j,k,\ell}|^p m_{j,k,\ell}^p\, \mu_{j,k,\ell}^{1-p} \right)^{1/p}.
\notag
\end{align}
with obvious modification for $p=\infty$.

Let us finally specialize the previous theorem to $\H_\A = \Co L^2_\A$.
In order to have a (Hilbert) frame in the usual sense, we have to renormalize the frame
elements, see also \cite[Theorem 4.7.7]{Rau_Diss}.

\begin{Theorem} Assume $g \in (S_0)_{rad} = (M^{1,1}_0)_{rad}$, $g\neq 0$. 
Then there exist constants $a_0,b_0>0$ such that for all $a, b$ satisfying $0<a<a_0,0<b<b_0$ the functions
\[
\tilde{g}_{j,k,\ell} \,:=\, \sqrt{\mu_{j,k,\ell}} ~\Omega\left(aj,bk, \sin \frac{\pi \ell}{2 N(j,k)}\right) g,
\qquad (j,k,\ell) \in I
\]
form a (Hilbert-) frame for $L^2_{rad}(\R^d)$.
\end{Theorem}

\subsection{Embeddings of radial modulation spaces}

Let us now apply the abstract results about embeddings of coorbit spaces to
modulation spaces of radial distributions. For simplicity we specialize to the spaces
$(M^p_s)_{rad} = (M^p_{m_s})_{rad}$, $1 \leq p \leq \infty$, $s \in \R$, 
with the weight function $m_s$ defined in (\ref{def_ms}).

\begin{Theorem}\label{thm_embed_mod} Let $1\leq p \leq q \leq \infty$ and $s,t \in \R$.
Set $\alpha:=1/p-1/q > 0$. 
\begin{itemize} 
\item[(a)] We have the continuous embedding $(M^p_s)_{rad}(\R^d) \hookrightarrow (M^q_t)_{rad}(\R^d)$ 
if and only if 
\[
t-s \,\leq\, \alpha(d-1).
\]
\item[(a)] The embedding $(M^p_s)_{rad}(\R^d) \hookrightarrow (M^q_t)_{rad}(\R^d)$ is compact 
if and only if 
\[
p<q \qquad \mbox{and } \qquad t-s \,<\, \alpha (d-1).
\]
\end{itemize}
\end{Theorem}
\begin{Proof} According to Theorem \ref{thm_compact} we need to 
investigate the sequence
\[
h(j,k,\ell) =
\frac{m_t(x_{j,k,\ell})}{m_s(x_{j,k,\ell})|V_{a,b}(x_{j,k,\ell},\omega_{j,k,\ell})|^{1/p-1/q}}
\,=\, (1+bk)^{t-s}|V_{a,b}(x_{j,k,\ell},\omega_{j,k,\ell})|^{-\alpha}.
\] 
Hence, by Theorem~\ref{thm_measVab} we get the estimation
\[
h(j,k,\ell) \,\leq\, C(1+j)^{-(d-1)\alpha}(1+k)^{t-s-(d-1)\alpha}, \quad (j,k,\ell) \in I
\]
Thus, if $t-s-(d-1)\alpha\leq 0$ then $h$ is contained in $\ell^\infty$, and if we have the
strict inequality $t-s-(d-1)\alpha < 0$ and $p<q$ then $h \in c_0$. This shows the
"if"-part of (a) and (b). On the other hand we have
\[
|V_{a,b}(x_{0,k,N(0,k)},\omega_{0,k,N(0,k)})| \, \asymp (1+k)^{d-1}.
\] 
Thus, if $t-s-(d-1)\alpha> 0$ then $h \notin \ell^\infty$ and if 
$t-s-(d-1)\alpha \geq 0$ or $p=q$ then $h \notin c_0$. This shows the "only if"-part.
\end{Proof}
We note the interesting special cases
\begin{align}
(M^1_{-(d-1)/2})_{rad}(\R^d) \,&\hookrightarrow\, (M^2_0)_{rad}(\R^d) \,=\, L^2_{rad}(\R^d),\notag\\
(S_0)_{rad}(\R^d) \,=\, (M^1_0)_{rad}(\R^d) \,&\hookrightarrow\, (M^2_{(d-1)/2})_{rad}(\R^d) \,=\, H^{(d-1)/2}_{rad}(\R^d),
\end{align} 
where $H^s(\R^d)$ denotes the Bessel potential space (Sobolev space) of index $s$. 
Moreover, $(S_0)_{rad}(\R^d)$ is compactly embedded into $L^2_{rad}(\R^d)$ for
$d \geq 2$.

In particular, the previous theorem shows that for $d\geq 2$, $s \in \R$ and $p < q$ the embedding
$(M^p_s)_{rad}(\R^d) \hookrightarrow (M^q_s)_{rad}(\R^d)$ is compact although
$M^p_s(\R^d)$ is not compactly embedded into $M^q_s(\R^d)$. Moreover, if $t$ is such that
$0<(t-s)<(d-1)(1/p-1/q)$ then $(M^p_s)_{rad}$ is (compactly) embedded into $(M^q_t)_{rad}$ although
$M^p_s$ is not even embedded into $M^q_t$.  

So roughly speaking, symmetry enforces
compactness of embeddings or even generates embeddings.
The first phenomenon was also observed for the Besov and Triebel-Lizorkin spaces
\cite{KLSS,SickS,LS,ST}
while the second phenomenon does not seem to be noticed earlier.

Also Theorem \ref{thm_pgq} can be applied to the modulation spaces. However, since
it is not specific to radial functions we omit its application here. We only note that
one cannot work with the weight function $m_s$ (because of the integrability condition). 
One rather has to take the functions $v_s(x,\omega) = (1+|x|+|\omega|)^s$ punishing
also the space variable.

Let us now gain more information about the compact embedding 
$(M^p_s)_{rad} \hookrightarrow (M^q_s)_{rad}$, $p < q$. Indeed,
let us compute the entropy numbers of the embedding operator. For an operator 
$S \in \B(E,F)$ the $n$-th entropy number $e_n(S)$ is defined by \cite{Pietsch}
\[
e_n(S) \,:=\, \inf
\{\epsilon > 0:\, \exists\, y_1,\hdots, y_{2^n-1} \in F \mbox{ such that } 
S(U_E) \subset \cup_{j=1}^{2^{n}-1} y_i + \epsilon U_F\}
\]
where $U_E, U_F$ denote the unit balls in $E,F$.
Moreover, we also need Lorentz sequence spaces. 
For $0<p \leq \infty$ the Lorentz space $\ell_{p,\infty}$ 
(sometimes also called weak $\ell^p$)
consists of those bounded sequences $\lambda=(\lambda_i)$ such that the quasi-norm
\[
\|\lambda\|_{p,\infty} \,:=\, \sup_{n\in \N} n^{1/p} s_n(\lambda)
\]
is finite, where $s_n(\lambda)$ denotes the non-increasing rearrangement of $\lambda$ defined
in (\ref{def_incre}).
It is well known that the class $\L^{(e)}_{p,\infty}$ of all operators $S$ 
whose sequence $(e_n(S))_{n \in \N}$ of entropy numbers is contained in $\ell_{p,\infty}$
forms an operator ideal \cite[Chapter 14.3]{Pietsch}, \cite[p.~35]{CS}.

The following theorem is a special case of Proposition 2 in \cite{Carl}.

\begin{Theorem}\label{thm_entropy} Let $1\leq p \leq q \leq \infty$ and $0<r\leq \infty$ be given and set
$1/s = 1/r + 1/p - 1/q$. Then $(w_i)_{i\in I} \in \ell_{r,\infty}$ implies that the embedding
operator $\Id : \ell^p \to \ell^q_{w}$ is contained in $\L^{(e)}_{s,\infty}$, i.e.,
$e_n(\Id : \ell^p \to \ell^q_w) \leq C n^{-1/s}$.
\end{Theorem}

In order to apply this theorem, we first need an estimate of the non-increasing rearrangement of 
the sequence $(|V_{a,b}(x_{j,k,\ell},\omega_{j,k,\ell})|^{-1})_{(j,k,\ell) \in I}$. 

\begin{Lemma}\label{lem_sn} Let $b_{j,k,\ell} := |V_{a,b}(x_{j,k,\ell},\omega_{j,k,\ell})|^{-1}$, $(j,k,\ell) \in I$.
Then it holds
\[
s_n(b) \,\leq\, C n^{-\frac{d-1}{3}},
\] 
i.e., $b \in \ell_{3/(d-1),\infty}$.
\end{Lemma}
\begin{Proof} Let us first count the number of indices in $I_n:=\{(j,k,\ell) \in I, j+k \leq n\}$ for $n \in \N_0$.
By (\ref{Njk_asymp}) we get
\begin{align}
\# I_n &\,=\, \sum_{j,k\geq 0, j+k \leq n} 2N(j,k) + 1 
\,\asymp\, 2n-1 + \sum_{j,k \geq 1, j+k \leq n} \frac{jk}{j+k}\notag\\
&=\, 2n-1 + \sum_{\ell=1}^n \frac{1}{\ell} \sum_{j=1}^{\ell-1} j(\ell-j)
\,\asymp\, n^3.\notag
\end{align}
Moreover, it follows from Theorem \ref{thm_measVab} that
\[
\sup_{(j,k,\ell) \in I \setminus I_n} b(j,k,\ell) \,\leq\, C n^{-(d-1)}
\] 
and, hence, $s_{\# I_n+1}(b) \leq C n^{-(d-1)}$. Let $\kappa$ denote the inverse function of $n \mapsto \# I_n$. 
Since $s_n(b)$ is a non-increasing sequence we obtain
\[
s_n(b) \,\leq\, C \kappa(n)^{-(d-1)} \leq C'n^{-\frac{d-1}{3}}.
\]
This completes the proof.  
\end{Proof}

Now we are ready to give an estimate of the entropy numbers of embeddings of radial
modulation spaces.

\begin{Theorem} Let $1\leq p < q \leq \infty$, $d\geq 2$ and $m$ be some moderate
invariant weight function
on $\R^d \times \R^d$. Then it holds
\[
e_n(\Id: (M^p_m)_{rad}(\R^d) \to (M^q_m)_{rad}(\R^d) ) \,\leq \, C n^{-\frac{d+2}{3}(1/p-1/q)}, \quad n \in \N.
\]
\end{Theorem}
\begin{Proof} Since $\L^{(e)}_{s,\infty}$ is an operator ideal, it suffices to estimate
the entropy numbers of the embedding $J$ from $(L^p_{m})^\bdiscrete_{SO(d)} = \ell^p_{m_p}$ into
$(L^q_{m})^\bdiscrete_{SO(d)} = \ell^q_{m_q}$ with 
$m_p(i) \,=\, m(z_i) |V_{a,b}(z_i)|^{1/p}$ by Theorem \ref{thm_coorbit_embed}.
Clearly, $J$ can be factorized
as 
\begin{equation}\label{factor_J}
J: \ell^p_{m_p} \to \ell^p \to \ell^{q}_v \to \ell^q_{m_q},\quad 
J \,=\, D_\sigma^{-1} J' D_{\sigma}
\end{equation}
where $D_\sigma$ is the (formal) diagonal operator $(D_\sigma x)_i = \sigma_i x_i$ with 
$\sigma_i = m_p(i)$. Furthermore,
\[
v(i) \,=\, \frac{m_q(i)}{m_p(i)} \,=\, |V_{a,b}(z_i)|^{-(1/p-1/q)}
\]
and $J' = \Id : \ell^p(I) \to \ell^q_v(I)$. 
By Lemma \ref{lem_sn} we have $v \in \ell_{r,\infty}$ with $1/r=\frac{d-1}{3}(1/p-1/q)$.
Theorem \ref{thm_entropy} yields $J' \in \L^{(e)}_{s,\infty}$ with
\[
1/s \,=\, 1/r + 1/p-1/q \,=\, \left(\frac{d-1}{3}+1\right)(1/p-1/q) \,=\, \frac{d+2}{3}(1/p-1/q).
\]
Since $\L^{(e)}_{s,\infty}$ is an operator ideal also $J \in \L^{(e)}_{s,\infty}$ by the factorization
(\ref{factor_J}). This concludes the proof.
\end{Proof}

\begin{Remark} Another measure of compactness of operators is provided by the approximation numbers.
For an operator $S \in \L(E,F)$ they are defined by
\[
a_n(S) \,=\, \inf \{\|S - T\|, T \in \L(E,F), \rank(T) \leq n\}.
\]
It is easy to see
that for weight functions $v,w$ on $I$ and $1\leq p\leq q \leq \infty$ 
it holds
\[
a_n(\Id: \ell^p_w \to \ell^q_v) \,\leq\, s_{n+1}\left(\left(\frac{v(i)}{w(i)}\right)_{i\in I}\right). 
\] 
Moreover, the class of operators $S$ whose sequence of approximation numbers $(a_n(S))_{n\in \N}$ is contained
in $\ell_{s,\infty}$ forms again an operator ideal. So with similar arguments as in the previous proof
one can show that
\[
a_n(\Id : (M^p_m)_{rad}(\R^d) \to (M^q_m)_{rad}(\R^d)) \,\leq\, C n^{-\frac{d-1}{3}(1/p-1/q)}
\]
for $1\leq p < q \leq \infty$, $d\geq 2$ and some moderate weight function $m$.
\end{Remark}

\subsection{Linear and nonlinear approximation}

Let us finally apply the abstract results on linear and nonlinear approximation to
the radial Gabor-like atomic decompositions of Theorem \ref{rad_mod_frame}. 
Let us first state the theorem for linear approximation (recall also the definition 
(\ref{def_linapprox_error}) of the error of linear approximation).

\begin{Theorem} Let $1\leq p < q \leq \infty$, $d\geq 2$ and $m$ be some moderate weight function.
Moreover, let $a,b> 0$ and $g$ such that 
$(\Omega(aj,bk,\sin \frac{\pi \ell}{2N(j,k)}) g)_{(j,k,\ell) \in I} = (g_{j,k,\ell})$
forms an atomic decomposition of $(M^p_m)_{rad}(\R^d)$ and of $(M^q_m)_{rad}$. 
Further, let $\tau:\N \to I$ be an ordering that realizes the non-increasing 
rearrangement of 
$(h_{j,k,\ell})=(|V_{a,b}(x_{j,k,\ell},\omega_{j,k,\ell})|^{-1}),~{(j,k,\ell) \in I}$,
i.e., $h_{\tau(n)} = s_n(h)$. Denote 
$V_n := \spann \{\Omega(x_{\tau(j)},\omega_{\tau(j)})g,
j=1,\hdots,n\}$. Then
\[
e(f,V_n,(M^q_m)_{rad}) \,\leq\, C n^{-\frac{d-1}{3}(1/p-1/q)} \|f|(M^p_m)_{rad}\|
\]
for all $f \in (M^p_m)_{rad}(\R^d)$.
\end{Theorem}
\begin{Proof} The claim follows immediately from Theorem \ref{thm_linapprox} together
with Lemma \ref{lem_sn}.
\end{Proof}

Of course, an important special case is $q=2$ and $m=1$ which corresponds to approximation
in $(M^2)_{rad} = L^2_{rad}(\R^d)$.
The theorem shows, in particular, that it is advantageous to approximate a radial function
with radial Gabor frames rather than with usual Gabor frames, 
see also Remark \ref{rem_lin}.
 
Let us now consider non-linear approximation. We assume once more 
that $a,b>0$ and $g$ are chosen such that
\[
(g_{j,k,\ell})_{(j,k,\ell)\in I} \,=\, (\Omega(aj,bk,\sin \frac{\pi \ell}{2N(j,k)}) g)_{(j,k,\ell) \in I} 
\]
forms an atomic decomposition for any of the spaces that we consider. Further we denote
\[
\sigma_n(f,(M^q_t)_{rad}) \,:=\, \inf_{N \subset I, \#N \leq n} 
\|f - \sum_{(j,k,\ell) \in N} \lambda_{j,k,\ell} g_{j,k,\ell} |(M^q_t)_{rad}\|
\]
the error of best $n$-term approximation in $M^q_t$, $1\leq q \leq \infty$, $t \in \R$.

\begin{Theorem} Let $1\leq p < q \leq \infty$ and $s,t \in \R$ and set $\alpha:= 1/p-1/q > 0$ and
assume 
\[
s \,\geq\, t- (d-1)(1/p-1/q) \,=\, t - \alpha(d-1),
\]
i.e., $(M^p_s)_{rad}(\R^d) \hookrightarrow (M^q_t)_{rad}(\R^d)$ by Theorem \ref{thm_embed_mod}.
Then
\begin{equation}\label{nonlin_radGab}
\left(\sum_{n=1}(n^\alpha \sigma_n(f, (M^q_t)_{rad}))^p\right)^{1/p} \leq C \|f|(M^p_s)_{rad}\|
\quad \mbox{for all }f \in (M^p_s)_{rad}.
\end{equation}
\end{Theorem}
\begin{Proof} This follows from Theorem \ref{thm_nonlinapprox} in connection with
Theorem \ref{thm_embed_mod}.
\end{Proof}

We remark that in the corresponding result for approximating with usual Gabor frames the condition
on $t,s$ would be $s\geq t$. However, for $d \geq 2$ we are allowed to choose $s = t - \alpha (d-1) < t$
when approximating with the radial Gabor frames. Thus, the class $(M^p_s)_{rad}$ for which 
(\ref{nonlin_radGab}) is guaranteed to hold for a certain $\alpha$ is larger than the one
for approximating an invariant $f$ with usual Gabor frames. Let us illustrate this for the special
case $p=1, q= 2$ and $t=0$ corresponding to approximation in $L^2_{rad}(\R^d)$. From the previous theorem
it follows that 
\[
\sigma_n(f, L^2_{rad}(\R^d)) \leq C \|f|(M^1_{-(d-1)/2})_{rad}\|\, n^{-1/2} \quad \mbox{for all }f \in
(M^1_{-(d-1)/2})_{rad}(\R^d).
\]
Now denote the error of best $n$-term approximation with usual Gabor frames 
by
\[
\sigma^*_n(f,M^{q}_t) \,=\, \inf_{N \subset \Z^d\times \Z^d, \# N \leq n} 
\|f - \sum_{(j,k) \in N } \lambda_{j,k} M_{bk} T_{aj} g|M^q_t\|     
\]
where the constants $a,b>0$ and the function $g$ are chosen according to Theorem \ref{thm_ab}.
We conclude from Theorem \ref{thm_nonlinapprox} (choosing $\A$ as the trivial automorphism group)
that
\[
\sigma^*_n(f,L^2(\R^d)) \,\leq\, C \|f|M^1\|\, n^{-1/2} \quad \mbox{for all } f\in M^1(\R^d) = S_0(\R^d),
\]   
in particular for all $f \in M^1_{rad}(\R^d)$. This is the best result we can get
from Theorem \ref{thm_nonlinapprox} (see also \cite{GS}) and there is no reason
why it should extent to the larger space $(M^1_{-(d-1)/2})_{rad}(\R^d))$. 

This shows to some extent that the radial Gabor frames perform better than 
usual Gabor frames when approximating radial functions.

\section{Acknowledgements}

The author would like to thank Prof. H.G. Feichtinger for pointing out
the papers \cite{KLSS,SickS,LS,ST}, which motivated to consider also
compactness of embeddings of radial modulation spaces. 
This paper was partly written during a stay at the Mathematical Institute
of the University of Wroc{\l}aw. He would like to thank its members for
the warm hospitality. The stay was supported by the European Union's Human 
Potential Programme under contract HPRN--CT--2001--00273 (HARP). The author also thanks 
the graduate program ``Applied Algorithmic Mathematics'' at the 
Technical University of Munich
(funded by the DFG) for its support.

\end{document}